\documentclass[graybox]{svmult}

\usepackage{type1cm}        
\usepackage{makeidx}         
\usepackage{graphicx}        
\usepackage{multicol}        
\usepackage[bottom]{footmisc}

\usepackage{newtxtext}       %
\usepackage{newtxmath}       

\makeindex             

\usepackage{amsmath}
\usepackage{tikz}
\usetikzlibrary{decorations.pathreplacing}
\usetikzlibrary{fit}
\newcommand\addvmargin[1]{
  \node[fit=(current bounding box),inner ysep=#1,inner xsep=0]{};
}
\usepackage{pgfplots}

\begin{document}

\bibliographystyle{plain}

\title*{On two subclasses of Motzkin paths and their relation to ternary trees}
\author{Helmut Prodinger \and Sarah J. Selkirk \and Stephan Wagner}
\authorrunning{H. Prodinger \and S. Selkirk \and S. Wagner}
\institute{Helmut Prodinger \and Sarah Selkirk \and Stephan Wagner \at Stellenbosch University, Department of Mathematical Sciences, Stellenbosch, 7602, \email{$\{\text{hproding, sjselkirk, swagner}\}$@sun.ac.za}}

%
%
\motto{To Peter Paule, a vibrant mathematician and a true innovator, on the occasion of his 60th birthday.} 
\maketitle

\abstract{Two subclasses of Motzkin paths, S-Motzkin and T-Motzkin paths, are introduced. We provide bijections between S-Motzkin paths and ternary trees, \mbox{S-Motzkin} paths and non-crossing trees, and T-Motzkin paths and ordered pairs of ternary trees. Symbolic equations for both paths, and thus generating functions for the paths, are provided. Using these, various parameters involving the two paths are analyzed.}

\section{Introduction}
\label{sec:1}

\vspace{1mm}\noindent
A Motzkin path is a non-negative lattice path with steps from the step set $\{\,
\begin{tikzpicture}[scale = 0.25, line width = 0.2mm]
\coordinate (aux) at (0,0);
\foreach \i in {0}
    \draw[line cap = round] (aux)--++(1,\i) coordinate (aux);
\end{tikzpicture}\,, \,
\begin{tikzpicture}[scale = 0.25, line width = 0.2mm]
\coordinate (aux) at (0,0);
\foreach \i in {1}
    \draw[line cap = round] (aux)--++(1,\i) coordinate (aux);
\end{tikzpicture}\,, \,
\begin{tikzpicture}[scale = 0.25, line width = 0.2mm]
\coordinate (aux) at (0,0);
\foreach \i in {-1}
    \draw[line cap = round] (aux)--++(1,\i) coordinate (aux);
\end{tikzpicture}\,
\}$ such that the path starts and ends on the $x$-axis. By placing further restrictions on Motzkin paths we obtain an interesting subclass.
\begin{definition}
\label{def:S-Motzkin}
An \emph{S-Motzkin path} is a Motzkin path of length $3n$ with $n$ of each type of step such that the following conditions hold
\begin{enumerate}
\item The initial step must be 
\begin{tikzpicture}[scale = 0.25, line width = 0.2mm]
\coordinate (aux) at (0,0);
\foreach \i in {0}
    \draw[line cap = round] (aux)--++(1,\i) coordinate (aux);
\end{tikzpicture}\,, and 
\item \begin{tikzpicture}[scale = 0.25, line width = 0.2mm]
\coordinate (aux) at (0,0);
\foreach \i in {0}
    \draw[line cap = round] (aux)--++(1,\i) coordinate (aux);
\end{tikzpicture} 
and  
\begin{tikzpicture}[scale = 0.25, line width = 0.2mm]
\coordinate (aux) at (0,0);
\foreach \i in {1}
    \draw[line cap = round] (aux)--++(1,\i) coordinate (aux);
\end{tikzpicture} steps alternate.
\end{enumerate}
\end{definition}
This definition was inspired by a question at the recent International Mathematics Competition \cite{IMC} involving restricted three-dimensional walks which can be translated into the two-dimensional S-Motzkin paths. These paths are enumerated by the generalized Catalan number, $\frac{1}{2n+1}\binom{3n}{n}$, and thus are bijective to ternary trees and non-crossing trees, as well as many other combinatorial objects \cite{Deutsch, Gu, Panholz, Prod, OEIS}. \mbox{We define} another subclass of Motzkin path which is related to both S-Motzkin paths and ternary trees.
\begin{definition}
A \emph{T-Motzkin path} is a Motzkin path of length $3n$ with $n$ of each type of step such that
\begin{enumerate}
\item The initial step is\, 
\begin{tikzpicture}[scale = 0.25, line width = 0.2mm]
\coordinate (aux) at (0,0);
\foreach \i in {1}
    \draw[line cap = round] (aux)--++(1,\i) coordinate (aux);
\end{tikzpicture}\,, and
\item 
\begin{tikzpicture}[scale = 0.25, line width = 0.2mm]
\coordinate (aux) at (0,0);
\foreach \i in {1}
    \draw[line cap = round] (aux)--++(1,\i) coordinate (aux);
\end{tikzpicture}\,
 and \,
\begin{tikzpicture}[scale = 0.25, line width = 0.2mm]
\coordinate (aux) at (0,0);
\foreach \i in {0}
    \draw[line cap = round] (aux)--++(1,\i) coordinate (aux);
\end{tikzpicture}\, steps
 alternate.
\end{enumerate}
\end{definition}
Note that although similar in definition, the class of T-Motzkin paths is larger than the class of S-Motzkin paths. Interchanging the\, 
\begin{tikzpicture}[scale = 0.25, line width = 0.2mm]
\coordinate (aux) at (0,0);
\foreach \i in {0}
    \draw[line cap = round] (aux)--++(1,\i) coordinate (aux);
\end{tikzpicture}\,
and\, 
\begin{tikzpicture}[scale = 0.25, line width = 0.2mm]
\coordinate (aux) at (0,0);
\foreach \i in {1}
    \draw[line cap = round] (aux)--++(1,\i) coordinate (aux);
\end{tikzpicture} 
steps in an arbitrary S-Motzkin path provides a T-Motzkin path, but the converse is not true. 
T-Motzkin paths of length $3n$ are enumerated by $\frac{1}{n+1}\binom{3n+1}{n}$ and thus bijective to the class of ordered pairs of ternary trees introduced by Knuth \cite{Knuth_Christmas}. There are several other equinumerous objects which can be found on the Online Encyclopedia of Integer Sequences A006013 \cite{OEIS}. 

Introducing another type of path is necessary for finding generating function equations for S-Motzkin and T-Motzkin paths, and thus we define a \emph{U-path} to be an S-Motzkin path without the initial\, 
\begin{tikzpicture}[scale = 0.25, line width = 0.2mm]
\coordinate (aux) at (0,0);
\foreach \i in {0}
    \draw[line cap = round] (aux)--++(1,\i) coordinate (aux);
\end{tikzpicture}\, step. Symbolic equations for T-Motzkin paths and U-paths can be obtained in terms of each other by making use of a decomposition based on the first return of the path. Since S-Motzkin paths and U-paths are `almost' the same, the generating function for S-Motzkin paths can be easily obtained from that of U-paths. 

Various parameters associated with different types of lattice paths have been studied \cite{Deutsch2, Pinzani, Wagner} and we provide analysis of the number of returns, peaks, valleys, and valleys on the $x$-axis in both S-Motzkin and T-Motzkin paths. This analysis is done using the symbolic equations and generating functions that are derived, as well as methods from the seminal book \emph{Analytic Combinatorics} by Flajolet and Sedgewick \cite{Flajolet}. During this analysis some interesting identities were found and are discussed briefly in Section \ref{sec:6}.

The study of these paths as well as parameters related to them has resulted in some generalizations and developments which will be reported in further publications. 

\section{Bijections}
\label{sec:2}

\vspace{1mm}\noindent

\subsection{S-Motzkin paths and ternary trees}
\label{sec:SMtoTT}

A bijection between S-Motzkin paths of length $3n$ and ternary trees with $n$ nodes is provided. 

\subsubsection{S-Motzkin paths to ternary trees}

We define $\varnothing$ to be the empty path. For an arbitrary S-Motzkin path $\mathcal{M}$, the canonical decomposition is given by
\begin{equation*}
\Phi(\mathcal{M}) = \big( \mathcal{A},\, \mathcal{B},\, \mathcal{C}\big),
\end{equation*}
where $\mathcal{A}$, $\mathcal{B}$, and $\mathcal{C}$ represent the S-Motzkin paths associated with the left, middle, and right subtrees respectively. Furthermore,  
\begin{itemize}
\item $\mathcal{C}$ is the path from the penultimate to the final return of $\mathcal{M}$, with the initial and final step removed,
\item $\mathcal{A}$ is the path from $y$ to $x$ (not including $x$), where $x$ is the first 
\begin{tikzpicture}[scale = 0.25, line width = 0.2mm]
\coordinate (aux) at (0,0);
\foreach \i in {0}
    \draw[line cap = round] (aux)--++(1,\i) coordinate (aux);
\end{tikzpicture} 
to the left of $\mathcal{C}$, $y$ is a  
\begin{tikzpicture}[scale = 0.25, line width = 0.2mm]
\coordinate (aux) at (0,0);
\foreach \i in {0}
    \draw[line cap = round] (aux)--++(1,\i) coordinate (aux);
\end{tikzpicture} step,  and the path from y to x is a Motzkin path of maximal length, and
\item $\mathcal{B}$ is what remains of $\mathcal{M}$ after removing the path from the penultimate to the final return of $\mathcal{M}$, as well as the path from $y$ to $x$ (including $x$).
\end{itemize}

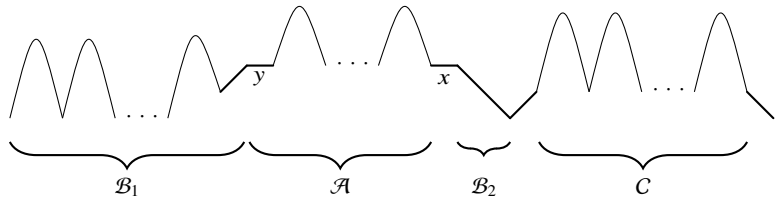
\begin{figure}[h]
\label{canonical}
	\begin{center}
		\begin{tikzpicture}[scale = 0.35]		
		\draw (2-2,0) .. controls (3-2,4) .. (4-2,0);
		\draw (4-2,0) .. controls (5-2,4) .. (6-2,0);
		\node at (5,0) {$\cdot$};
		\node at (5.5,0) {$\cdot$};
		\node at (4.5,0) {$\cdot$};
		\draw (6,0) .. controls (7,4) .. (8,1);
		\draw [thick, line cap = round](8,1) to  (9,2);
		\draw [thick, line cap = round](9,2) to  (10,2);
		
		\node at (5.5+7,2) {$\cdot$};
		\node at (6+7,2) {$\cdot$};
		\node at (6.5+7,2) {$\cdot$};
		\draw (10,2) .. controls (11,5) .. (12,2);	
		\draw (14,2) .. controls (15,5) .. (16,2);	

		\draw [thick, line cap = round](16,2) to  (17,2);
		\draw [thick, line cap = round](17,2) to  (19,0);
		\draw [thick, line cap = round](19,0) to  (20,1);
		\draw (19+1,1) .. controls (20+1,5) .. (21+1,1);
		\draw (21+1,1) .. controls (22+1,5) .. (23+1,1);
		\node at (23.5+1,1) {$\cdot$};
		\node at (24+1,1) {$\cdot$};
		\node at (24.5+1,1) {$\cdot$};
		\draw (25+1,1) .. controls (26+1,5) .. (27+1,1);
		\draw [thick, line cap = round](28,1) to  (29,0);
		
		\node at (9.5,1.5) {$y$};
		\node at (16.5,1.5) {$x$};

		\draw [thick, decorate, decoration={brace, amplitude=10pt, mirror, raise=4pt}] (0cm, -0.5) to
		node[below,yshift=-0.5cm] {$\mathcal{B}_{1}$}
		(8.9cm, -0.5);
		
		\draw [thick, decorate, decoration={brace, amplitude=5pt, mirror, raise=4pt}] (17cm, -0.5) to
		node[below,yshift=-0.5cm] {$\mathcal{B}_{2}$}
		(19cm, -0.5);

		\draw [thick, decorate, decoration={brace, amplitude=10pt, mirror, raise=4pt}] (9.1cm, -0.5) to
		node[below,yshift=-0.5cm] {$\mathcal{A}$}
		(16cm,- 0.5);
		
		\draw [thick, decorate, decoration={brace, amplitude=10pt, mirror, raise=4pt}] (20.1cm, -0.5) to
		node[below,yshift=-0.5cm] {$\mathcal{C}$}
		(28cm, -0.5);
		
		\end{tikzpicture}
		
	\end{center}
	\caption{Canonical decomposition of an arbitrary S-Motzkin path}

\end{figure}

This process is performed recursively and terminates at an empty path. Note that each application of $\Phi$ adds one node and removes one of each type of step. This proves inductively that an S-Motzkin path of length $3n$ maps to a ternary tree with $n$ (internal) nodes.

\subsubsection{Ternary trees to S-Motzkin paths}

The inverse mapping is performed recursively on the end nodes as follows. Each node of a ternary tree has three (possibly empty) subtrees. Call the paths associated with the left, middle, and right subtrees $\mathcal{A}$, $\mathcal{B}$, and $\mathcal{C}$ respectively. 

Starting at the end nodes, replace each node with $
\mathcal{B}_{1}\,\mathcal{A}\,
\begin{tikzpicture}[scale = 0.25, line width = 0.2mm]
\coordinate (aux) at (0,0);
\foreach \i in {0}
    \draw[line cap = round] (aux)--++(1,\i) coordinate (aux);
\end{tikzpicture}\,
\mathcal{B}_{2}\,
\begin{tikzpicture}[scale = 0.25, line width = 0.2mm]
\coordinate (aux) at (0,0);
\foreach \i in {1}
    \draw[line cap = round] (aux)--++(1,\i) coordinate (aux);
\end{tikzpicture}\,
\mathcal{C} \,
\begin{tikzpicture}[scale = 0.25, line width = 0.2mm]
\coordinate (aux) at (0,0);
\foreach \i in {-1}
    \draw[line cap = round] (aux)--++(1,\i) coordinate (aux);
\end{tikzpicture}\,,
$ 
where $\mathcal{B}_{1}$ is the path from the start of $\mathcal{B}$ to the final 
\begin{tikzpicture}[scale = 0.25, line width = 0.2mm]
\coordinate (aux) at (0,0);
\foreach \i in {1}
    \draw[line cap = round] (aux)--++(1,\i) coordinate (aux);
\end{tikzpicture}
step of $\mathcal{B}$. The path $\mathcal{B}_{2}$ is what remains of $\mathcal{B}$ after removing $\mathcal{B}_{1}$.
This process is continued recursively on each set of end nodes and terminates at the root to produce an \mbox{S-Motzkin} path. Note that for each node that is removed one of each type of step is added, and thus a ternary tree with $n$ nodes produces an S-Motzkin path of length $3n$. 

\subsubsection{Example}

As an example, we map the following S-Motzkin path into a ternary tree. Since the steps are reversible, the inverse mapping can be seen by reading the example in reverse.
Let $\mathcal{M}$ be the S-Motzkin path
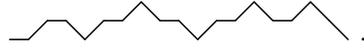
\begin{center}
\begin{tikzpicture}[scale = 0.25, line width = 0.2mm]
\coordinate (aux) at (0,0);
\foreach \i in {0, 1, 0, -1, 1, 0, 1, -1, 0, -1, 1, 0, 1, -1, 0, 1, -1, -1}
    \draw[line cap = round] (aux)--++(1,\i) coordinate (aux);
\end{tikzpicture}\, .
\end{center}
The canonical decomposition of $\mathcal{M}$ is then 
$\Phi(\mathcal{M}) = \Big( 
\begin{tikzpicture}[scale = 0.2, line width = 0.2mm]
\coordinate (aux) at (0,0);
\foreach \i in {0, 1, -1}
    \draw[line cap = round] (aux)--++(1,\i) coordinate (aux);
\end{tikzpicture}
\,, \, 
\begin{tikzpicture}[scale = 0.2, line width = 0.2mm]
\coordinate (aux) at (0,0);
\foreach \i in {0, 1, 0, -1, 1, -1}
    \draw[line cap = round] (aux)--++(1,\i) coordinate (aux);
\end{tikzpicture}	
\,, \,
\begin{tikzpicture}[scale = 0.2, line width = 0.2mm]
\coordinate (aux) at (0,0);
\foreach \i in {0, 1, -1, 0, 1, -1}
    \draw[line cap = round] (aux)--++(1,\i) coordinate (aux);
\end{tikzpicture}
\Big).$
Hence
\begin{center}
\scalebox{0.75}{
\begin{tikzpicture}
\node (1) at (0, 0){$\bullet$};

\node (2) at (-2, -1){
\begin{tikzpicture}[scale = 0.2, line width = 0.2mm]
\coordinate (aux) at (0,0);
\foreach \i in {0, 1, -1}
    \draw[line cap = round] (aux)--++(1,\i) coordinate (aux);
\end{tikzpicture}
};
\node (3) at (0, -1){
\begin{tikzpicture}[scale = 0.2, line width = 0.2mm]
\coordinate (aux) at (0,0);
\foreach \i in {0, 1, 0, -1, 1, -1}
    \draw[line cap = round] (aux)--++(1,\i) coordinate (aux);
\end{tikzpicture}	
};
\node(4) at (2, -1){
\begin{tikzpicture}[scale = 0.2, line width = 0.2mm]
\coordinate (aux) at (0,0);
\foreach \i in {0, 1, -1, 0, 1, -1}
    \draw[line cap = round] (aux)--++(1,\i) coordinate (aux);
\end{tikzpicture}
};

\draw[-] (1.center) to (2);
\draw[-] (1.center) to (3);
\draw[-] (1.center) to (4);
\end{tikzpicture}}
\end{center}
Continuing recursively:

\begin{minipage}{0.3\textwidth}
\centering
\scalebox{0.7}{
\begin{tikzpicture}
\node (1) at (0, 0){$\bullet$};
\node (2) at (-1.5, -1){$\bullet$};
\node (3) at (0, -1){$\bullet$};
\node (4) at (1.5, -1){$\bullet$};
\node (5) at (-2, -2){$\varnothing$};
\node (6) at (-1.5, -2){$\varnothing$};
\node (7) at (-1, -2){$\varnothing$};
\node (8) at (-0.5, -2){$\varnothing$};
\node (9) at (0, -2){
\begin{tikzpicture}[scale = 0.2, line width = 0.2mm]
\coordinate (aux) at (0,0);
\foreach \i in {0, 1, -1}
    \draw[line cap = round] (aux)--++(1,\i) coordinate (aux);
\end{tikzpicture}
};
\node (10) at (0.5, -2){$\varnothing$};
\node (11) at (1, -2){
\begin{tikzpicture}[scale = 0.2, line width = 0.2mm]
\coordinate (aux) at (0,0);
\foreach \i in {0, 1, -1}
    \draw[line cap = round] (aux)--++(1,\i) coordinate (aux);
\end{tikzpicture}
};
\node (12) at (1.5, -2){$\varnothing$};
\node (13) at (2, -2){$\varnothing$};

\draw[-] (1.center) to (2.center);
\draw[-] (1.center) to (3.center);
\draw[-] (1.center) to (4.center);
\draw[-] (2.center) to (5);
\draw[-] (2.center) to (6);
\draw[-] (2.center) to (7);
\draw[-] (3.center) to (8);
\draw[-] (3.center) to (9);
\draw[-] (3.center) to (10);
\draw[-] (4.center) to (11);
\draw[-] (4.center) to (12);
\draw[-] (4.center) to (13);
\end{tikzpicture}}
\end{minipage}\quad$\rightarrow$
\begin{minipage}{.3\textwidth}
\centering
\scalebox{0.7}{
\begin{tikzpicture}
\node (1) at (0, 0){$\bullet$};
\node (2) at (-1.3, -1){$\bullet$};
\node (3) at (0, -1){$\bullet$};
\node (4) at (1.3, -1){$\bullet$};
\node (8) at (0, -1.6){$\bullet$};
\node (11) at (1, -1.6){$\bullet$};

\node (14) at (0.3, -2.5){$\varnothing$};
\node (15) at (0, -2.5){$\varnothing$};
\node (16) at (-0.3, -2.5){$\varnothing$};
\node (17) at (0.7, -2.5){$\varnothing$};
\node (18) at (1, -2.5){$\varnothing$};
\node (19) at (1.3, -2.5){$\varnothing$};

\draw[-] (1.center) to (2.center);
\draw[-] (1.center) to (3.center);
\draw[-] (1.center) to (4.center);
\draw[-] (3.center) to (8.center);
\draw[-] (4.center) to (11.center);

\draw[-] (8.center) to (14);
\draw[-] (8.center) to (15);
\draw[-] (8.center) to (16);
\draw[-] (11.center) to (17);
\draw[-] (11.center) to (18);
\draw[-] (11.center) to (19);
\end{tikzpicture}}
\end{minipage}
$\rightarrow$
\begin{minipage}{.3\textwidth}
\centering
\scalebox{0.7}{
\begin{tikzpicture}
\node (1) at (0, 0){$\bullet$};
\node (2) at (-1.3, -1){$\bullet$};
\node (3) at (0, -1){$\bullet$};
\node (4) at (1.3, -1){$\bullet$};
\node (8) at (0, -1.6){$\bullet$};
\node (11) at (1, -1.6){$\bullet$};

\draw[-] (1.center) to (2.center);
\draw[-] (1.center) to (3.center);
\draw[-] (1.center) to (4.center);
\draw[-] (3.center) to (8.center);
\draw[-] (4.center) to (11.center);
\end{tikzpicture}}
\end{minipage}


\begin{table}[h!]
\caption{Bijection for $n=3$}
\setlength{\tabcolsep}{5mm}
\def\arraystretch{1.25}
\centering
\begin{tabular}{|c|c||c|c|}
  \hline
  S-Motzkin path   &   Ternary tree & S-Motzkin path   &   Ternary tree 
  \\ \hline
\begin{tikzpicture}[scale = 0.25, line width = 0.2mm]
\coordinate (aux) at (0,0);
\foreach \i in {0, 1, -1, 0, 1, -1, 0, 1, -1}
    \draw[line cap = round] (aux)--++(1,\i) coordinate (aux);
\addvmargin{1mm}
\end{tikzpicture}  &    
\begin{tikzpicture}
\node (1) at (0, 0){$\bullet$};
\node (2) at (-1/2, -1/2){$\bullet$};
\node (5) at (-2/2, -2/2){$\bullet$};
\draw[-] (1.center) to (2.center);
\draw[-] (2.center) to (5.center);
\end{tikzpicture}
&


\begin{tikzpicture}[scale = 0.25, line width = 0.2mm]
\coordinate (aux) at (0,0);
\foreach \i in {0, 1, 0, 1, -1, 0, 1, -1, -1}
    \draw[line cap = round] (aux)--++(1,\i) coordinate (aux);
\addvmargin{1mm}
\end{tikzpicture}  & 
\begin{tikzpicture}
\node (1) at (0, 0){$\bullet$};
\node (4) at (1/2, -1/2){$\bullet$};
\node (7) at (0, -2/2){$\bullet$};

\draw[-] (4.center) to (7.center);
\draw[-] (1.center) to (4.center);
\end{tikzpicture}  
  \\ \hline

  \begin{tikzpicture}[scale = 0.25, line width = 0.2mm]
\coordinate (aux) at (0,0);
\foreach \i in {0, 1, -1, 0, 1, 0, -1, 1, -1}
    \draw[line cap = round] (aux)--++(1,\i) coordinate (aux);
\addvmargin{1mm}
\end{tikzpicture}
   &   
\begin{tikzpicture}
\node (1) at (0, 0){$\bullet$};
\node (3) at (0, -1/2){$\bullet$};
\node (6) at (-1/2, -2/2){$\bullet$};

\draw[-] (1.center) to (3.center);
\draw[-] (3.center) to (6.center);
\end{tikzpicture}
&

\begin{tikzpicture}[scale = 0.25, line width = 0.2mm]
\coordinate (aux) at (0,0);
\foreach \i in {0, 1, 0, 1, -1, 0, -1, 1, -1}
    \draw[line cap = round] (aux)--++(1,\i) coordinate (aux);
\addvmargin{1mm}
\end{tikzpicture}   & 
\begin{tikzpicture}
\node (1) at (0, 0){$\bullet$};
\node (2) at (-1/2, -1/2){$\bullet$};
\node (3) at (0, -1/2){$\bullet$};
\draw[-] (1.center) to (2.center);
\draw[-] (1.center) to (3.center);
\end{tikzpicture}
   
  \\ \hline

\begin{tikzpicture}[scale = 0.25, line width = 0.2mm]
\coordinate (aux) at (0,0);
\foreach \i in {0, 1, -1, 0, 1, 0, 1, -1, -1}
    \draw[line cap = round] (aux)--++(1,\i) coordinate (aux);
\addvmargin{1mm}
\end{tikzpicture}   &  
\begin{tikzpicture}
\node (1) at (0, 0){$\bullet$};
\node (2) at (-1/2, -1/2){$\bullet$};
\node (4) at (1/2, -1/2){$\bullet$};

\draw[-] (1.center) to (2.center);
\draw[-] (1.center) to (4.center);
\end{tikzpicture}
&

\begin{tikzpicture}[scale = 0.25, line width = 0.2mm]
\coordinate (aux) at (0,0);
\foreach \i in {0, 1, 0, 1, -1, -1, 0, 1, -1}
    \draw[line cap = round] (aux)--++(1,\i) coordinate (aux);
\addvmargin{1mm}
\end{tikzpicture}   & 
\begin{tikzpicture}
\node (1) at (0, 0){$\bullet$};
\node (2) at (-1/2, -1/2){$\bullet$};
\node (7) at (0, -2/2){$\bullet$};

\draw[-] (1.center) to (2.center);
\draw[-] (2.center) to (7.center);
\end{tikzpicture}

  \\ \hline

\begin{tikzpicture}[scale = 0.25, line width = 0.2mm]
\coordinate (aux) at (0,0);
\foreach \i in {0, 1, 0, -1, 1, -1, 0, 1, -1}
    \draw[line cap = round] (aux)--++(1,\i) coordinate (aux);
\addvmargin{1mm}
\end{tikzpicture}
   &  
\begin{tikzpicture}
\node (1) at (0, 0){$\bullet$};
\node (2) at (-1/2, -1/2){$\bullet$};
\node (6) at (-1/2, -2/2){$\bullet$};

\draw[-] (1.center) to (2.center);
\draw[-] (2.center) to (6.center);
\end{tikzpicture}
&

\begin{tikzpicture}[scale = 0.25, line width = 0.2mm]
\coordinate (aux) at (0,0);
\foreach \i in {0, 1, 0, 1, 0, -1, -1, 1, -1}
    \draw[line cap = round] (aux)--++(1,\i) coordinate (aux);
\addvmargin{1mm}
\end{tikzpicture}   & 
\begin{tikzpicture}
\node (1) at (0, 0){$\bullet$};
\node (3) at (0, -1/2){$\bullet$};
\node (8) at (1/2, -2/2){$\bullet$};

\draw[-] (1.center) to (3.center);
\draw[-] (3.center) to (8.center);
\end{tikzpicture}

  \\ \hline

\begin{tikzpicture}[scale = 0.25, line width = 0.2mm]
\coordinate (aux) at (0,0);
\foreach \i in {0, 1, 0, -1, 1, 0, -1, 1, -1}
    \draw[line cap = round] (aux)--++(1,\i) coordinate (aux);
\addvmargin{1mm}
\end{tikzpicture}
   &  
\begin{tikzpicture}
\node (1) at (0, 0){$\bullet$};
\node (3) at (0, -1/2){$\bullet$};
\node (7) at (0, -2/2){$\bullet$};

\draw[-] (1.center) to (3.center);
\draw[-] (3.center) to (7.center);
\end{tikzpicture}
&

\begin{tikzpicture}[scale = 0.25, line width = 0.2mm]
\coordinate (aux) at (0,0);
\foreach \i in {0, 1, 0, 1, 0, -1, 1, -1, -1}
    \draw[line cap = round] (aux)--++(1,\i) coordinate (aux);
\addvmargin{1mm}
\end{tikzpicture}  &  
\begin{tikzpicture}
\node (1) at (0, 0){$\bullet$};
\node (4) at (1/2, -1/2){$\bullet$};
\node (8) at (1/2, -2/2){$\bullet$};

\draw[-] (1.center) to (4.center);
\draw[-] (4.center) to (8.center);
\end{tikzpicture}

  \\ \hline

\begin{tikzpicture}[scale = 0.25, line width = 0.2mm]
\coordinate (aux) at (0,0);
\foreach \i in {0, 1, 0, -1, 1, 0, 1, -1, -1}
    \draw[line cap = round] (aux)--++(1,\i) coordinate (aux);
\addvmargin{1mm}
\end{tikzpicture}
   &  
\begin{tikzpicture}
\node (1) at (0, 0){$\bullet$};
\node (3) at (0, -1/2){$\bullet$};
\node (4) at (1/2, -1/2){$\bullet$};

\draw[-] (1.center) to (3.center);
\draw[-] (1.center) to (4.center);
\end{tikzpicture}
& 

\begin{tikzpicture}[scale = 0.25, line width = 0.2mm]
\coordinate (aux) at (0,0);
\foreach \i in {0, 1, 0, 1, 0, 1, -1, -1, -1}
    \draw[line cap = round] (aux)--++(1,\i) coordinate (aux);
\addvmargin{1mm}
\end{tikzpicture}  &  
\begin{tikzpicture}
\node (1) at (0, 0){$\bullet$};
\node (4) at (1/2, -1/2){$\bullet$};
\node (9) at (2/2, -2/2){$\bullet$};

\draw[-] (4.center) to (9.center);
\draw[-] (1.center) to (4.center);
\end{tikzpicture}
\\ \hline

  \end{tabular}
\end{table}

\subsection{T-Motzkin paths and pairs of ternary trees}
\label{sec:bijection_pairs}

\subsubsection{T-Motzkin paths to pairs of ternary trees}
Since a bijection between S-Motzkin paths and ternary trees is already provided, we show that every T-Motzkin path can be decomposed uniquely into an ordered pair of S-Motzkin paths (possibly including an empty path). 

Given an arbitrary T-Motzkin path $\mathcal{N}$, we perform a canonical decomposition $\Omega(\mathcal{N}) = (\mathcal{A}, \mathcal{B})$ where
\begin{itemize}
\item $\mathcal{B}$ is the path from $y$ to $x$ (not including $x$) where $x$ is the rightmost 
\begin{tikzpicture}[scale = 0.25, line width = 0.2mm]
\coordinate (aux) at (0,0);
\foreach \i in {0}
    \draw[line cap = round] (aux)--++(1,\i) coordinate (aux);
\end{tikzpicture} 
step of $\mathcal{N}$, $y$ is a  
\begin{tikzpicture}[scale = 0.25, line width = 0.2mm]
\coordinate (aux) at (0,0);
\foreach \i in {0}
    \draw[line cap = round] (aux)--++(1,\i) coordinate (aux);
\end{tikzpicture}
step, and the path from $y$ to $x$ is a Motzkin path of maximal length, and
\item $\mathcal{A}$ is what remains of $\mathcal{N}$ after removing the path from $y$ to $x$ (including $x$), with an additional 
\begin{tikzpicture}[scale = 0.25, line width = 0.2mm]
\coordinate (aux) at (0,0);
\foreach \i in {0}
    \draw[line cap = round] (aux)--++(1,\i) coordinate (aux);
\end{tikzpicture}
step at the start of the path. In Figure \ref{canonicalT} this is the path 
\begin{tikzpicture}[scale = 0.25, line width = 0.2mm]
\coordinate (aux) at (0,0);
\foreach \i in {0}
    \draw[line cap = round] (aux)--++(1,\i) coordinate (aux);
\end{tikzpicture}\,$\mathcal{A}_{1}\mathcal{A}_{2}$.
\end{itemize}
Note that both $\mathcal{A}$ and $\mathcal{B}$ are S-Motzkin paths.
\begin{figure}[h]
	\begin{center}
		\begin{tikzpicture}[scale = 0.3]
		\draw [thick, line cap = round](0, 0) to (1, 1);

		\draw (1,1) .. controls (2,5) .. (3,1);
		\draw (3,1) .. controls (4,5) .. (5,1);
		\node at (5.5,1) {$\cdot$};
		\node at (6,1.4) {$\cdot$};
		\node at (6.5,1.8) {$\cdot$};
		\draw (7,2) .. controls (8,6) .. (9,2);

		\draw [thick, line cap = round](9, 2) to (10, 3);
		\draw [thick, line cap = round](10, 3) to (11, 3);
		\draw [thick, line cap = round](11, 3) to (12, 4);

		\draw (12,4) .. controls (13,8) .. (14,4);
		\draw (14,4) .. controls (15,8) .. (16,4);
		\node at (16.5,4) {$\cdot$};
		\node at (17,4) {$\cdot$};
		\node at (17.5,4) {$\cdot$};
		\draw (18,4) .. controls (19,8) .. (20,4);
		
		\draw [thick, line cap = round](20, 4) to (21, 3);
		\draw [thick, line cap = round](21, 3) to (22, 3);
		\draw [thick, line cap = round](22, 3) to (23, 2);
		\node at (23.2,1.7) {$\cdot$};
		\node at (23.5,1.4) {$\cdot$};
		\node at (23.8, 1.1) {$\cdot$};
		\draw [thick, line cap = round](24, 1) to (25, 0);

		\draw [thick, decorate, decoration={brace, amplitude=5pt, mirror, raise=4pt}] (0cm, -0.5) to
		node[below,yshift=-0.5cm] {$\mathcal{A}_{1}$}(9.9cm, -0.5);

		\draw [thick, decorate, decoration={brace, amplitude=5pt, mirror, raise=4pt}] (10cm, -0.5) to
		node[below,yshift=-0.5cm] {$\mathcal{B}$}(21cm, -0.5);

		\draw [thick, decorate, decoration={brace, amplitude=5pt, mirror, raise=4pt}] (22.1cm, -0.5) to
		node[below,yshift=-0.5cm] {$\mathcal{A}_{2}$}(25cm, -0.5);

		\node at (21.5, 2.5){$x$};
		\node at (10.5, 2.5){$y$};
		\end{tikzpicture}
		
	\end{center}
	\caption{Canonical decomposition of an arbitrary T-Motzkin path}

\label{canonicalT}
\end{figure}

\subsubsection{Pairs of ternary trees to T-Motzkin paths}

Given an arbitrary pair of ternary trees, we can use the bijection given in Section \ref{sec:SMtoTT} to obtain an ordered pair of S-Motzkin paths, $(\mathcal{A}, \mathcal{B})$. All S-Motzkin paths start with a 
\begin{tikzpicture}[scale = 0.25, line width = 0.2mm]
\coordinate (aux) at (0,0);
\foreach \i in {0}
    \draw[line cap = round] (aux)--++(1,\i) coordinate (aux);
\end{tikzpicture}
step and 
end in an 
\begin{tikzpicture}[scale = 0.25, line width = 0.2mm]
\coordinate (aux) at (0,0);
\foreach \i in {1}
    \draw[line cap = round] (aux)--++(1,\i) coordinate (aux);
\end{tikzpicture}
step followed by a series of\, 
\begin{tikzpicture}[scale = 0.25, line width = 0.2mm]
\coordinate (aux) at (0,0);
\foreach \i in {-1}
    \draw[line cap = round] (aux)--++(1,\i) coordinate (aux);
\end{tikzpicture}\,
steps. To obtain a T-Motzkin path from $(\mathcal{A}, \mathcal{B})$ we
\begin{itemize}
\item remove the initial 
\begin{tikzpicture}[scale = 0.25, line width = 0.2mm]
\coordinate (aux) at (0,0);
\foreach \i in {0}
    \draw[line cap = round] (aux)--++(1,\i) coordinate (aux);
\end{tikzpicture}
step from $\mathcal{A}$, and
\item insert the path $\mathcal{B}\begin{tikzpicture}[scale = 0.25, line width = 0.2mm]
\coordinate (aux) at (0,0);
\foreach \i in {0}
    \draw[line cap = round] (aux)--++(1,\i) coordinate (aux);
\end{tikzpicture}$
immediately after the final 
\begin{tikzpicture}[scale = 0.25, line width = 0.2mm]
\coordinate (aux) at (0,0);
\foreach \i in {1}
    \draw[line cap = round] (aux)--++(1,\i) coordinate (aux);
\end{tikzpicture}
step of $\mathcal{A}$.
\end{itemize}	
\subsubsection{Example}

We provide an example of the mapping from T-Motzkin paths to ternary trees. The inverse mapping can be seen by reading this example in reverse. 
Let $\mathcal{N}$ be 
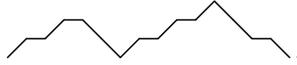
\begin{equation*}
\begin{tikzpicture}[scale = 0.25, line width = 0.2mm]
\coordinate (aux) at (0,0);
\foreach \i in {1, 0, 1, 0, -1, -1, 1, 0, 1, 0, 1, -1, -1, 0, -1	}
    \draw[line cap = round] (aux)--++(1,\i) coordinate (aux);
\end{tikzpicture}\,.
\end{equation*}
Then $\Omega(\mathcal{N})$ is given by 
\begin{equation*}
\raisebox{2ex}{$\Biggl($}
\raisebox{1ex}{
\begin{tikzpicture}[scale = 0.25, line width = 0.2mm]
\coordinate (aux) at (0,0);
\foreach \i in {0, 1, 0, 1, 0, -1, -1, 1, -1}
    \draw[line cap = round] (aux)--++(1,\i) coordinate (aux);
\end{tikzpicture}\, ,\,  
\begin{tikzpicture}[scale = 0.25, line width = 0.2mm]
\coordinate (aux) at (0,0);
\foreach \i in {0, 1, 0, 1, -1, -1}
    \draw[line cap = round] (aux)--++(1,\i) coordinate (aux);
\end{tikzpicture}}
\raisebox{2ex}{$\Biggr)$} \qquad\qquad \raisebox{2.5ex}{$\longrightarrow$} \qquad\qquad
\raisebox{2ex}{$\Biggl($}
\scalebox{0.6}{
\begin{tikzpicture}
\node (1) at (0, 0){$\bullet$};
\node (2) at (1/2, -1){$\bullet$};
\node (3) at (0, -1/2){$\bullet$};
\draw[-] (3.center) to (2.center);
\draw[-] (1.center) to (3.center);
\end{tikzpicture}}\, , \,
\scalebox{0.6}{
\begin{tikzpicture}
\node (1) at (0, 0){$\bullet$};
\node (2) at (1/2, -1/2){$\bullet$};
\draw[-] (1.center) to (2.center);
\end{tikzpicture}}
\raisebox{2ex}{$\Biggr)$}.
\end{equation*}

\subsection{S-Motzkin paths and non-crossing trees}

The definition of a non-crossing tree as well as the representation of a non-crossing tree that is used in this text (with a marker to separate left and right children) can be found in \cite{Panholz}. To assist in describing the bijection, we define a \emph{piece} to be a maximal subpath of a Motzkin path consisting of (in order) one up step, a series of down steps (possibly empty), one horizontal step, and a series of down steps (possibly empty). Note that an arbitrary S-Motzkin path of length $3n$ consists of an initial 
\begin{tikzpicture}[scale = 0.25, line width = 0.2mm]
\coordinate (aux) at (0,0);
\foreach \i in {0}
    \draw[line cap = round] (aux)--++(1,\i) coordinate (aux);
\end{tikzpicture}\, step followed by $n-1$ pieces, and a final 
\begin{tikzpicture}[scale = 0.25, line width = 0.2mm]
\coordinate (aux) at (0,0);
\foreach \i in {1}
    \draw[line cap = round] (aux)--++(1,\i) coordinate (aux);
\end{tikzpicture}\, step followed by a series of 
\begin{tikzpicture}[scale = 0.25, line width = 0.2mm]
\coordinate (aux) at (0,0);
\foreach \i in {-1}
    \draw[line cap = round] (aux)--++(1,\i) coordinate (aux);
\end{tikzpicture}\, steps.
Each piece is uniquely determined by the number of\, 
\begin{tikzpicture}[scale = 0.25, line width = 0.2mm]
\coordinate (aux) at (0,0);
\foreach \i in {-1}
    \draw[line cap = round] (aux)--++(1,\i) coordinate (aux);
\end{tikzpicture}\, steps and the position of the 
\begin{tikzpicture}[scale = 0.25, line width = 0.2mm]
\coordinate (aux) at (0,0);
\foreach \i in {0}
    \draw[line cap = round] (aux)--++(1,\i) coordinate (aux);
\end{tikzpicture}\, step. The \emph{characteristic pair} of a piece is the ordered pair $(t, i)$ with $t$ denoting the number of\, 
\begin{tikzpicture}[scale = 0.25, line width = 0.2mm]
\coordinate (aux) at (0,0);
\foreach \i in {-1}
    \draw[line cap = round] (aux)--++(1,\i) coordinate (aux);
\end{tikzpicture}\, steps in the piece, and $i$ denoting the position of the 
\begin{tikzpicture}[scale = 0.25, line width = 0.2mm]
\coordinate (aux) at (0,0);
\foreach \i in {0}
    \draw[line cap = round] (aux)--++(1,\i) coordinate (aux);
\end{tikzpicture}\, step (with the\, 
\begin{tikzpicture}[scale = 0.25, line width = 0.2mm]
\coordinate (aux) at (0,0);
\foreach \i in {1}
    \draw[line cap = round] (aux)--++(1,\i) coordinate (aux);
\end{tikzpicture}\, step in position $0$).
\begin{figure}[h]
\begin{center}
\begin{tikzpicture}[scale = 0.25, line width = 0.2mm]
\coordinate (aux) at (0,0);
\foreach \i in {0, 1, -1, 0, 1, 0, 1, -1, 0, -1, 1, 0, 1, -1, -1}
    \draw[line cap = round] (aux)--++(1,\i) coordinate (aux);
\draw [decorate, decoration={brace, amplitude=2pt, mirror, raise=4pt}] (1cm, -0.5) to node[below,yshift=-0.25cm] {$1$}(4cm, -0.5);
\draw [decorate, decoration={brace, amplitude=2pt, mirror, raise=4pt}] (4.1cm, -0.5) to node[below,yshift=-0.25cm] {$2$}(6cm, -0.5);
\draw [decorate, decoration={brace, amplitude=2pt, mirror, raise=4pt}] (6.1cm, -0.5) to node[below,yshift=-0.25cm] {$3$}(10cm, -0.5);
\draw [decorate, decoration={brace, amplitude=2pt, mirror, raise=4pt}] (10.1cm, -0.5) to node[below,yshift=-0.25cm] {$4$}(12cm, -0.5);
\end{tikzpicture}
\end{center}
\caption{The four pieces in the given S-Motzkin path of length $15$}
\label{fig:pieces}
\end{figure}
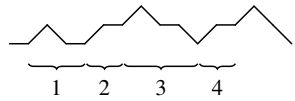

\subsubsection{S-Motzkin paths to non-crossing trees}

Given an arbitrary S-Motzkin path of length $3n$, we let $c$ denote the number of steps in the final series of\, 
\begin{tikzpicture}[scale = 0.25, line width = 0.2mm]
\coordinate (aux) at (0,0);
\foreach \i in {-1}
    \draw[line cap = round] (aux)--++(1,\i) coordinate (aux);
\end{tikzpicture}\, steps of the path. In the resulting non-crossing tree, draw a root with $c$ children.
Considering the pieces of the path from right to left, we find the characteristic pair $(t, i)$ and draw $i-1$ left subtrees and $t-i+1$ right subtrees on the rightmost available node.

\subsubsection{Non-crossing trees to S-Motzkin paths}

Given an arbitrary non-crossing tree, let $c$ denote the number of children of the root. Associate with each non-root node an ordered pair $(u, j+1)$ where $u$ equals the number of children of the node and $j$ equals the number of left subtrees of the node. Then remove the leftmost leaf. Draw an initial 
\begin{tikzpicture}[scale = 0.25, line width = 0.2mm]
\coordinate (aux) at (0,0);
\foreach \i in {0}
    \draw[line cap = round] (aux)--++(1,\i) coordinate (aux);
\end{tikzpicture}\, step and repeat the following until only the root remains: consider the leftmost leaf's ordered pair $(u, j+1)$ and add an 
\begin{tikzpicture}[scale = 0.25, line width = 0.2mm]
\coordinate (aux) at (0,0);
\foreach \i in {1}
    \draw[line cap = round] (aux)--++(1,\i) coordinate (aux);
\end{tikzpicture}\, step, $j$ 
\begin{tikzpicture}[scale = 0.25, line width = 0.2mm]
\coordinate (aux) at (0,0);
\foreach \i in {-1}
    \draw[line cap = round] (aux)--++(1,\i) coordinate (aux);
\end{tikzpicture}\, steps, a 
\begin{tikzpicture}[scale = 0.25, line width = 0.2mm]
\coordinate (aux) at (0,0);
\foreach \i in {0}
    \draw[line cap = round] (aux)--++(1,\i) coordinate (aux);
\end{tikzpicture}\, step and $u-j$ 
\begin{tikzpicture}[scale = 0.25, line width = 0.2mm]
\coordinate (aux) at (0,0);
\foreach \i in {-1}
    \draw[line cap = round] (aux)--++(1,\i) coordinate (aux);
\end{tikzpicture}\, steps to the path, and then remove the leftmost leaf from the tree. Finally, add an 
\begin{tikzpicture}[scale = 0.25, line width = 0.2mm]
\coordinate (aux) at (0,0);
\foreach \i in {1}
    \draw[line cap = round] (aux)--++(1,\i) coordinate (aux);
\end{tikzpicture}\, step, and $c$ 
\begin{tikzpicture}[scale = 0.25, line width = 0.2mm]
\coordinate (aux) at (0,0);
\foreach \i in {-1}
    \draw[line cap = round] (aux)--++(1,\i) coordinate (aux);
\end{tikzpicture}\, steps to the path.


\begin{table}[h!]
\caption{Bijection for $n=3$}
\setlength{\tabcolsep}{3mm}
\def\arraystretch{1.2}
\centering
\begin{tabular}{|c|c||c|c|}
  \hline
  S-Motzkin path   &   Non-crossing tree & S-Motzkin path   &   Non-crossing tree 
  \\ \hline
\begin{tikzpicture}[scale = 0.25, line width = 0.2mm]
\coordinate (aux) at (0,0);
\foreach \i in {0, 1, -1, 0, 1, -1, 0, 1, -1}
    \draw[line cap = round] (aux)--++(1,\i) coordinate (aux);
\addvmargin{1mm}
\end{tikzpicture}  &    
\scalebox{0.4}{
\begin{tikzpicture}[every node/.style={circle,inner sep=1.5pt,minimum size=1.5em}]
\node(1)[draw, circle,inner sep=1.5pt] at (0.0, 0.0){};
\node(2)[draw, circle,inner sep=1.5pt]  at (0.0, -1.0){};
\node(3)[draw, circle,inner sep=1.5pt]  at (-1.0,-2.0){};
\node(4)[draw, circle,inner sep=1.5pt]  at (-2.0,-3.0){};

\draw[-]   (1) to (2);
\draw[-]   (2) to (3);
\draw[-]   (3) to (4);

\draw[->,>=stealth,ultra thick](0,-1.3) to (2);
\draw[->,>=stealth,ultra thick](-1,-2.3) to (3);
\end{tikzpicture}}
&


\begin{tikzpicture}[scale = 0.25, line width = 0.2mm]
\coordinate (aux) at (0,0);
\foreach \i in {0, 1, 0, 1, -1, 0, 1, -1, -1}
    \draw[line cap = round] (aux)--++(1,\i) coordinate (aux);
\addvmargin{1mm}
\end{tikzpicture}  & 
\scalebox{0.4}{
\begin{tikzpicture}[every node/.style={circle,inner sep=1.5pt,minimum size=1.5em}]
\node(1)[draw, circle,inner sep=1.5pt] at (0.0, 0.0){};
\node(2)[draw, circle,inner sep=1.5pt]  at (1.0, -1.0){};
\node(3)[draw, circle,inner sep=1.5pt]  at (-1.0,-1.0){};
\node(4)[draw, circle,inner sep=1.5pt]  at (0.0,-2.0){};

\draw[-]   (1) to (2);
\draw[-]   (1) to (3);
\draw[-]   (2) to (4);

\draw[->,>=stealth,ultra thick](1.0,-1.3) to (2);
\end{tikzpicture}} 
  \\ \hline

  \begin{tikzpicture}[scale = 0.25, line width = 0.2mm]
\coordinate (aux) at (0,0);
\foreach \i in {0, 1, -1, 0, 1, 0, -1, 1, -1}
    \draw[line cap = round] (aux)--++(1,\i) coordinate (aux);
\addvmargin{1mm}
\end{tikzpicture}
   &   
\scalebox{0.4}{
\begin{tikzpicture}[every node/.style={circle,inner sep=1.5pt,minimum size=1.5em}]
\node(1)[draw, circle,inner sep=1.5pt] at (0.0, 0.0){};
\node(2)[draw, circle,inner sep=1.5pt]  at (0.0, -1.0){};
\node(3)[draw, circle,inner sep=1.5pt]  at (1.0,-2.0){};
\node(4)[draw, circle,inner sep=1.5pt]  at (0.0,-3.0){};

\draw[-]   (1) to (2);
\draw[-]   (2) to (3);
\draw[-]   (3) to (4);

\draw[->,>=stealth,ultra thick](0,-1.3) to (2);
\draw[->,>=stealth,ultra thick](1,-2.3) to (3);
\end{tikzpicture}} 
&

\begin{tikzpicture}[scale = 0.25, line width = 0.2mm]
\coordinate (aux) at (0,0);
\foreach \i in {0, 1, 0, 1, -1, 0, -1, 1, -1}
    \draw[line cap = round] (aux)--++(1,\i) coordinate (aux);
\addvmargin{1mm}
\end{tikzpicture}   & 
\scalebox{0.4}{
\begin{tikzpicture}[every node/.style={circle,inner sep=1.5pt,minimum size=1.5em}]
\node(1)[draw, circle,inner sep=1.5pt] at (0.0, 0.0){};
\node(2)[draw, circle,inner sep=1.5pt]  at (0.0, -1.0){};
\node(3)[draw, circle,inner sep=1.5pt]  at (-1.0,-2.0){};
\node(4)[draw, circle,inner sep=1.5pt]  at (1.0,-2.0){};

\draw[-]   (1) to (2);
\draw[-]   (2) to (3);
\draw[-]   (2) to (4);

\draw[->,>=stealth,ultra thick](0,-1.3) to (2);
\end{tikzpicture}} 
   
  \\ \hline

\begin{tikzpicture}[scale = 0.25, line width = 0.2mm]
\coordinate (aux) at (0,0);
\foreach \i in {0, 1, -1, 0, 1, 0, 1, -1, -1}
    \draw[line cap = round] (aux)--++(1,\i) coordinate (aux);
\addvmargin{1mm}
\end{tikzpicture}   &  
\scalebox{0.4}{
\begin{tikzpicture}[every node/.style={circle,inner sep=1.5pt,minimum size=1.5em}]
\node(1)[draw, circle,inner sep=1.5pt] at (0.0, 0.0){};
\node(2)[draw, circle,inner sep=1.5pt]  at (1.0, -1.0){};
\node(3)[draw, circle,inner sep=1.5pt]  at (-1.0,-1.0){};
\node(4)[draw, circle,inner sep=1.5pt]  at (-2.0,-2.0){};

\draw[-]   (1) to (2);
\draw[-]   (1) to (3);
\draw[-]   (3) to (4);

\draw[->,>=stealth,ultra thick](-1.0,-1.3) to (3);
\end{tikzpicture}} 

&

\begin{tikzpicture}[scale = 0.25, line width = 0.2mm]
\coordinate (aux) at (0,0);
\foreach \i in {0, 1, 0, 1, -1, -1, 0, 1, -1}
    \draw[line cap = round] (aux)--++(1,\i) coordinate (aux);
\addvmargin{1mm}
\end{tikzpicture}   & 
\scalebox{0.4}{
\begin{tikzpicture}[every node/.style={circle,inner sep=1.5pt,minimum size=1.5em}]
\node(1)[draw, circle,inner sep=1.5pt] at (0.0, 0.0){};
\node(2)[draw, circle,inner sep=1.5pt]  at (0.0, -1.0){};
\node(3)[draw, circle,inner sep=1.5pt]  at (-1.0,-2.0){};
\node(4)[draw, circle,inner sep=1.5pt]  at (-2.0,-2.0){};

\draw[-]   (1) to (2);
\draw[-]   (2) to (3);
\draw[-]   (2) to (4);

\draw[->,>=stealth,ultra thick](0,-1.3) to (2);
\end{tikzpicture}} 

  \\ \hline

\begin{tikzpicture}[scale = 0.25, line width = 0.2mm]
\coordinate (aux) at (0,0);
\foreach \i in {0, 1, 0, -1, 1, -1, 0, 1, -1}
    \draw[line cap = round] (aux)--++(1,\i) coordinate (aux);
\addvmargin{1mm}
\end{tikzpicture}
   &  
\scalebox{0.4}{
\begin{tikzpicture}[every node/.style={circle,inner sep=1.5pt,minimum size=1.5em}]
\node(1)[draw, circle,inner sep=1.5pt] at (0.0, 0.0){};
\node(2)[draw, circle,inner sep=1.5pt]  at (0.0, -1.0){};
\node(3)[draw, circle,inner sep=1.5pt]  at (-1.0,-2.0){};
\node(4)[draw, circle,inner sep=1.5pt]  at (0.0,-3.0){};

\draw[-]   (1) to (2);
\draw[-]   (2) to (3);
\draw[-]   (3) to (4);

\draw[->,>=stealth,ultra thick](0,-1.3) to (2);
\draw[->,>=stealth,ultra thick](-1,-2.3) to (3);
\end{tikzpicture}} 

&

\begin{tikzpicture}[scale = 0.25, line width = 0.2mm]
\coordinate (aux) at (0,0);
\foreach \i in {0, 1, 0, 1, 0, -1, -1, 1, -1}
    \draw[line cap = round] (aux)--++(1,\i) coordinate (aux);
\addvmargin{1mm}
\end{tikzpicture}   & 
\scalebox{0.4}{
\begin{tikzpicture}[every node/.style={circle,inner sep=1.5pt,minimum size=1.5em}]
\node(1)[draw, circle,inner sep=1.5pt] at (0.0, 0.0){};
\node(2)[draw, circle,inner sep=1.5pt]  at (0.0, -1.0){};
\node(3)[draw, circle,inner sep=1.5pt]  at (1.0,-2.0){};
\node(4)[draw, circle,inner sep=1.5pt]  at (2.0,-2.0){};

\draw[-]   (1) to (2);
\draw[-]   (2) to (3);
\draw[-]   (2) to (4);

\draw[->,>=stealth,ultra thick](0,-1.3) to (2);
\end{tikzpicture}} 

  \\ \hline

\begin{tikzpicture}[scale = 0.25, line width = 0.2mm]
\coordinate (aux) at (0,0);
\foreach \i in {0, 1, 0, -1, 1, 0, -1, 1, -1}
    \draw[line cap = round] (aux)--++(1,\i) coordinate (aux);
\addvmargin{1mm}
\end{tikzpicture}
   &  
\scalebox{0.4}{
\begin{tikzpicture}[every node/.style={circle,inner sep=1.5pt,minimum size=1.5em}]
\node(1)[draw, circle,inner sep=1.5pt] at (0.0, 0.0){};
\node(2)[draw, circle,inner sep=1.5pt]  at (0.0, -1.0){};
\node(3)[draw, circle,inner sep=1.5pt]  at (1.0,-2.0){};
\node(4)[draw, circle,inner sep=1.5pt]  at (2.0,-3.0){};

\draw[-]   (1) to (2);
\draw[-]   (2) to (3);
\draw[-]   (3) to (4);

\draw[->,>=stealth,ultra thick](0,-1.3) to (2);
\draw[->,>=stealth,ultra thick](1,-2.3) to (3);
\end{tikzpicture}} 
&

\begin{tikzpicture}[scale = 0.25, line width = 0.2mm]
\coordinate (aux) at (0,0);
\foreach \i in {0, 1, 0, 1, 0, -1, 1, -1, -1}
    \draw[line cap = round] (aux)--++(1,\i) coordinate (aux);
\addvmargin{1mm}
\end{tikzpicture}  &  
\scalebox{0.4}{
\begin{tikzpicture}[every node/.style={circle,inner sep=1.5pt,minimum size=1.5em}]
\node(1)[draw, circle,inner sep=1.5pt] at (0.0, 0.0){};
\node(2)[draw, circle,inner sep=1.5pt]  at (1.0, -1.0){};
\node(3)[draw, circle,inner sep=1.5pt]  at (-1.0,-1.0){};
\node(4)[draw, circle,inner sep=1.5pt]  at (2.0,-2.0){};

\draw[-]   (1) to (2);
\draw[-]   (1) to (3);
\draw[-]   (2) to (4);

\draw[->,>=stealth,ultra thick](1,-1.3) to (2);
\end{tikzpicture}} 

  \\ \hline

\begin{tikzpicture}[scale = 0.25, line width = 0.2mm]
\coordinate (aux) at (0,0);
\foreach \i in {0, 1, 0, -1, 1, 0, 1, -1, -1}
    \draw[line cap = round] (aux)--++(1,\i) coordinate (aux);
\addvmargin{1mm}
\end{tikzpicture}
   &  
\scalebox{0.4}{
\begin{tikzpicture}[every node/.style={circle,inner sep=1.5pt,minimum size=1.5em}]
\node(1)[draw, circle,inner sep=1.5pt] at (0.0, 0.0){};
\node(2)[draw, circle,inner sep=1.5pt]  at (1.0, -1.0){};
\node(3)[draw, circle,inner sep=1.5pt]  at (-1.0,-1.0){};
\node(4)[draw, circle,inner sep=1.5pt]  at (0.0,-2.0){};

\draw[-]   (1) to (2);
\draw[-]   (1) to (3);
\draw[-]   (3) to (4);

\draw[->,>=stealth,ultra thick](-1.0,-1.3) to (3);
\end{tikzpicture}} 
& 

\begin{tikzpicture}[scale = 0.25, line width = 0.2mm]
\coordinate (aux) at (0,0);
\foreach \i in {0, 1, 0, 1, 0, 1, -1, -1, -1}
    \draw[line cap = round] (aux)--++(1,\i) coordinate (aux);
\addvmargin{1mm}
\end{tikzpicture}  &  
\scalebox{0.4}{
\begin{tikzpicture}[every node/.style={circle,inner sep=1.5pt,minimum size=1.5em}]
\node(1)[draw, circle,inner sep=1.5pt] at (0.0, 0.0){};
\node(2)[draw, circle,inner sep=1.5pt]  at (0.0, -1.0){};
\node(3)[draw, circle,inner sep=1.5pt]  at (-1.0,-1.0){};
\node(4)[draw, circle,inner sep=1.5pt]  at (1.0,-1.0){};

\draw[-]   (1) to (2);
\draw[-]   (1) to (3);
\draw[-]   (1) to (4);

\end{tikzpicture}}
\\ \hline

  \end{tabular}
\end{table}

\section{Generating functions and related paths}
\label{sec:3}

\vspace{1mm}\noindent
Let $\mathcal{T}$ be the class of T-Motzkin paths, $\mathcal{U}$ be the class of U-paths, and 
\begin{equation*}
T(z) = \sum_{n\geq 0}t_{n}z^{n} \qquad \text{ and } \qquad U(z) = \sum_{n\geq 0}u_{n}z^{n}
\end{equation*}
be their respective generating functions, where $t_{n}$ and $u_{n}$ represent the number of paths of length $n$ in the given class.

We derive symbolic equations for the two types of paths based on a first return decomposition. Note that the only U-path of length less than five is given by\,
\begin{tikzpicture}[scale = 0.25, line width = 0.2mm]
\coordinate (aux) at (0,0);
\foreach \i in {1, -1}
    \draw[line cap = round] (aux)--++(1,\i) coordinate (aux);
\end{tikzpicture}\,.
Taking into account the first return of a U-path, it is clear that a U-path of length five or more can be decomposed as either
\begin{equation*}
\text{\raisebox{1ex}{$(a)$}}\quad
\begin{tikzpicture}[scale=0.5]
\draw [line cap = round](0,0) to  (0.5,0.5);
\draw [line cap = round](0.5,0.5) to  (1,0.5);
	\node at (1.6, 0.8) {$\mathcal{X}$};	
\draw [line cap = round](2.1,0.5) to  (2.6,0);
	\node at (3, 0.3) {$\mathcal{Y}$};
\end{tikzpicture}\qquad \text{\raisebox{0.5ex}{ or }} \qquad
\text{\raisebox{1ex}{$(b)$}}\quad
\begin{tikzpicture}[scale=0.5]
\draw [line cap = round](0,0) to  (0.5,0.5);
\draw [line cap = round](0.5,0.5) to  (1,0);
\draw [line cap = round](1,0) to  (1.5,0);
	\node at (2, 0.3) {$\mathcal{Z}$};
\end{tikzpicture}.
\end{equation*}

In $(a)$, $\mathcal{X}$ can either be a U-path or a T-Motzkin path. If $\mathcal{X}$ is a U-path, then $\mathcal{Y}$ is either empty or $\mathcal{Y}$ is a\, 
\begin{tikzpicture}[scale = 0.25, line width = 0.2mm]
\coordinate (aux) at (0,0);
\foreach \i in {0}
    \draw[line cap = round] (aux)--++(1,\i) coordinate (aux);
\end{tikzpicture} step followed by a U-path. If $\mathcal{X}$ is a T-Motzkin path, then $\mathcal{Y}$ is a U-path.
In $(b)$, $\mathcal{Z}$ has to be a U-path. This then results in the symbolic equation
\begin{equation*}
\mathcal{U} \quad = \quad
\begin{tikzpicture}[scale=0.45]
\draw [line cap = round](0,0) to  (0.5,0.5);
\draw [line cap = round](0.5,0.5) to  (1,0);
\end{tikzpicture}
\quad + \quad
\begin{tikzpicture}[scale=0.45]
\draw [line cap = round](0,0) to  (0.5,0.5);
\draw [line cap = round](0.5,0.5) to  (1,0.5);
	\node at (1.6, 0.8) {$\mathcal{T}$};	
\draw [line cap = round](1.8,0.5) to  (2.3,0);
	\node at (2.8, 0.3) {$\mathcal{U}$};
\end{tikzpicture}
\quad + \quad
\begin{tikzpicture}[scale=0.45]
\draw [line cap = round](0,0) to  (0.5,0.5);
\draw [line cap = round](0.5,0.5) to  (1,0);
\draw [line cap = round](1,0) to  (1.5,0);
	\node at (2, 0.3) {$\mathcal{U}$};
\end{tikzpicture}
\quad + \quad
\begin{tikzpicture}[scale=0.45]
\draw [line cap = round](0,0) to  (0.5,0.5);
\draw [line cap = round](0.5,0.5) to  (1,0.5);
	\node at (1.6, 0.8) {$\mathcal{U}$};	
\draw [line cap = round](2.1,0.5) to  (2.6,0);
\end{tikzpicture}
\quad+\quad
\begin{tikzpicture}[scale=0.45]
\draw [line cap = round](0,0) to  (0.5,0.5);
\draw [line cap = round](0.5,0.5) to  (1,0.5);
	\node at (1.6, 0.8) {$\mathcal{U}$};	
\draw [line cap = round](2.1,0.5) to  (2.6,0);
\draw [line cap = round](2.6,0) to  (3.1,0);
	\node at (3.6, 0.3) {$\mathcal{U}$};
\end{tikzpicture},
\end{equation*}
from which we obtain the equation 
\begin{equation}
\label{eqn:1}
U(z) = z^{2} + z^{3}T(z)U(z) + 2z^{3}U(z) + z^{4}U(z)^{2}.
\end{equation}
Again, any T-Motzkin path of length $3n$ is either empty or, considering the first return of the path, of the form $(a)$ or $(b)$ as given in the U-path case.

In $(a)$, $\mathcal{X}$ can either be a U-path or a T-Motzkin path. If $\mathcal{X}$ is a U-path, then an `extra' 
\begin{tikzpicture}[scale = 0.25, line width = 0.2mm]
\coordinate (aux) at (0,0);
\foreach \i in {0}
    \draw[line cap = round] (aux)--++(1,\i) coordinate (aux);
\end{tikzpicture}\,
 step needs to appear in $\mathcal{Y}$, and thus $\mathcal{Y}$ is given by a 
\begin{tikzpicture}[scale = 0.25, line width = 0.2mm]
\coordinate (aux) at (0,0);
\foreach \i in {0}
    \draw[line cap = round] (aux)--++(1,\i) coordinate (aux);
\end{tikzpicture}\, step followed by a T-Motzkin path.
If $\mathcal{X}$ is a T-Motzkin path, then $\mathcal{Y}$ is also a T-Motzkin path. With analogous reasoning we can see that for $(b)$ the only possibility for $\mathcal{Z}$ is a T-Motzkin path. Using this we obtain the symbolic equation
\begin{equation*}
\mathcal{T} \quad = \quad \varepsilon \quad + \quad
\begin{tikzpicture}[scale=0.5]
\draw [line cap = round](0,0) to  (0.5,0.5);
\draw [line cap = round](0.5,0.5) to  (1,0.5);
	\node at (1.6, 0.8) {$\mathcal{T}$};	
\draw [line cap = round](1.8,0.5) to  (2.3,0);
	\node at (2.8, 0.3) {$\mathcal{T}$};
\end{tikzpicture}
\quad + \quad
\begin{tikzpicture}[scale=0.5]
\draw [line cap = round](0,0) to  (0.5,0.5);
\draw [line cap = round](0.5,0.5) to  (1,0);	
\draw [line cap = round](1,0) to  (1.5,0);
	\node at (2, 0.3) {$\mathcal{T}$};
\end{tikzpicture}
\quad + \quad
\begin{tikzpicture}[scale=0.5]
\draw [line cap = round](0,0) to  (0.5,0.5);
\draw [line cap = round](0.5,0.5) to  (1,0.5);
	\node at (1.6, 0.8) {$\mathcal{U}$};	
\draw [line cap = round](2.1,0.5) to  (2.6,0);
\draw [line cap = round](2.6,0) to  (3.1,0);
	\node at (3.7, 0.3) {$\mathcal{T}$};
\end{tikzpicture}
\end{equation*}
which results in the equation
\begin{equation}
\label{eqn:2}
T(z) = 1 + z^{3}T(z)^{2} + z^{3}T(z) + z^{4}U(z)T(z).
\end{equation}
Solving the system of equations given by (\ref{eqn:1}) and (\ref{eqn:2}) yields
\begin{align*}
T(z) &= 1 + 2z^{3}T(z)^{2}-z^{6}T(z)^{3},\\
U(z) & = z^{2} + 3z^{3}U(z) + 3z^{4}U(z)^{2}+z^{5}U(z)^{3}
\end{align*}
which, with substitutions, is amenable to application of the Lagrange inversion formula. To demonstrate this, consider the equation
\begin{align*}
T(z) & = 1 + 2z^{3}T(z)^{2} - z^{6}T(z)^{3}.
\end{align*}
This can be factorised as $T(z)(1-z^{3}T(z))^{2} = 1$, and with substitutions $R = z^{3}T(z)$ and $x = z^{3}$ we find that $x = R(1-R)^{2}$. Therefore 
\begin{align*}
[z^{3n}]T(z) & = [x^{n+1}]R = \frac{1}{n+1}[w^{n}]\frac{1}{\left(1-w\right)^{2n+2}} = \frac{1}{n+1}\binom{3n+1}{n},
\end{align*}
which results in 
\begin{equation*}
T(z) = \sum_{n\geq 0}\frac{1}{n+1}\binom{3n+1}{n}z^{3n} \quad \text{ and similarly } \quad U(z) = \sum_{n\geq 1}\frac{1}{2n+1}\binom{3n}{n}z^{3n-1}.
\end{equation*}
Since U-paths are S-Motzkin paths without the initial horizontal step, the generating function for S-Motzkin paths is given by
$S(z) = \sum_{n\geq 1}\frac{1}{2n+1}\binom{3n}{n}z^{3n}$. 
Note that \mbox{$(1+S(z))^{2} = T(z)$}, which was pointed out by Knuth in his 2014 Christmas lecture \cite{Knuth_Christmas}. We have proved this by means of the bijection provided in Section \ref{sec:bijection_pairs}.

\section{Analysis of various parameters}
\label{sec:5}

\vspace{1mm}\noindent
In this section the analysis of the number of returns is done in detail, and results for the number of peaks, the number of valleys, and the number of valleys on the $x$-axis are done similarly. The study of these parameters in Dyck paths can be found in \cite{Deutsch2, Mansour}.
\subsection{The number of returns}
From the generating functions for U-paths and T-Motzkin paths along with the substitutions $x = z^{3}$ and $x = t(1-t)^{2}$, we obtain 
\begin{equation*}
T(z)=\frac{1}{(1-t)^2}\qquad \text{and} \qquad S(z)=\frac{t}{1-t}.
\end{equation*}
We introduce the variable $u$ to count the number of returns, and from the symbolic equations for U-paths and T-Motzkin paths we obtain the bivariate generating functions:
\begin{align*}
S(z, u)&= u^{2}z^{3} + uz^{3}S(z, u)T(z, 1) + u^{2}z^{3}S(z, u) + u^{2}z^{3}S(z, 1) + u^{2}z^{3}S(z, 1)S(z, u),\\
T(z, u)&= 1+uz^{3}T(z, 1)T(z, u)+u^2z^{3}T(z, u)+u^2z^{3}S(z, 1)T(z, u).
\end{align*}
Solving this system of equations we find that 
\begin{equation*}
S(z, u) = \frac{{(1 - t)} t u^{2}}{1-tu -tu^{2} + t^{2} u^{2}} \qquad \text{and} \qquad T(z, u) = \frac{1}{1-tu-tu^{2} + t^{2} u^{2}}.
\end{equation*}

\subsubsection{Mean and variance}
\label{sec:5.1.1}
For a bivariate generating function $K(z, u)$ with $u$ representing the parameter of interest, we obtain the mean and variance as follows. The mean is given by
\begin{align*}
K_{\text{ave}} = [z^{n}]\frac{\partial}{\partial u}K(z, u)\Big|_{u=1} \, \Big/\, [z^{n}]K(z, 1),
\end{align*} 
and the variance is
\begin{align*}
K_{\text{var}} = [z^{n}]\frac{\partial^{2}}{(\partial u)^{2}}K(z, u)\Big|_{u=1}\, \Big/\, [z^{n}]K(z, 1) + K_{\text{ave}} - \big(K_{\text{ave}}\big)^{2}.
\end{align*}
In the sections that follow some simplifications occur when calculating variances. These are discussed in more detail in Section \ref{sec:6}.

To determine the average number of returns we calculate the derivative of $S(z, u)$ and $T(z, u)$ with respect to $u$, 
\begin{equation*}
\frac{\partial }{\partial u}S(z, u)\Big|_{u=1}=\frac{{(2 - t)} t}{{(1 - t)}^{3}}\qquad \text{ and } \qquad \frac{\partial }{\partial u}T(z, u)\Big|_{u=1}=\frac{t(3-2t)}{ (1-t)^4}.
\end{equation*}
The total number of returns in all paths of length $3n$ is then obtained by extracting the coefficients of these expressions by means of Cauchy's integral formula. For S-Motzkin paths this results in 
\begin{align*}
[x^n]\frac{{(2 - t)} t}{{(1 - t)}^{3}}
& = \frac{1}{2\pi i}\oint\frac{1}{(t(1-t)^{2})^{n+1}}\cdot\frac{t(2-t)}{(1-t)^3}\cdot (1-t)(1-3t) \, dt\\
& = \frac{1}{2\pi i}\oint\frac{1}{t^{n}}\cdot\frac{2-7t+3t^{2}}{(1-t)^{2n+4}}\, dt
  = [t^{n-1}]\frac{2-7t+3t^{2}}{(1-t)^{2n+4}}\\
&  = 2\binom{3n+2}{n-1} -7\binom{3n+1}{n-2} +3\binom{3n}{n-3},
\end{align*}
and for T-Motzkin paths we obtain
\begin{align*}
[x^n]\frac{t(3-2t)}{ (1-t)^4} & = 3\binom{3n+3}{n-1}-11\binom{3n+2}{n-2}+6\binom{3n+1}{n-3}.
\end{align*}
Therefore in S-Motzkin paths the average number of returns for paths of length $3n$ is
\begin{align*}
\frac{2\binom{3n+2}{n-1} -7\binom{3n+1}{n-2} +3\binom{3n}{n-3}}{\frac{1}{2n+1}\binom{3n}{n}} = \frac{n(23n+17)}{2 {(2 n + 3)} {(n + 1)}} = \frac{23}{4} - \frac{81}{8n} + \mathcal{O}\Big(\frac{1}{n^{2}}\Big)
\end{align*}
and for T-Motzkin paths the average number of returns is
\begin{align*}
\frac{3\binom{3n+3}{n-1}-11\binom{3n+2}{n-2}+6\binom{3n+1}{n-3}}{\frac{1}{n+1}\binom{3n+1}{n}} = \frac{{(19n + 26)} n}{2 (2n+3)(n+2)} = \frac{19}{4} - \frac{81}{8n} + \mathcal{O}\Big(\frac{1}{n^{2}}\Big).
\end{align*}

To calculate the variance in the number of returns for paths of length $3n$, we find the second derivatives of $S(z, u)$ and $T(z, u)$ with respect to $u$:
\begin{equation*}
\frac{\partial^{2} }{(\partial u) ^{2}}S(z, u)\Big|_{u=1}=\frac{2t {(1+3t-4t^{2}+t^{3})}}{{(1 - t)}^{5}}
\end{equation*}
and 
\begin{equation*}
\frac{\partial^{2} }{(\partial u) ^{2}}T(z, u)\Big|_{u=1} = \frac{2t {(1 + 6t -9t^{2} + 3t^{3})}}{{(1 - t)}^{6}}. 
\end{equation*}
We again determine the coefficients using Cauchy's integral formula, 
\begin{align*}
[x^{n}]\frac{2t(1+3t-4t^{2}+t^{3})}{(1-t)^{5}} & = 2\bigg[\binom{3n+4}{n-1} - 13\binom{3n+2}{n-3} + 13\binom{3n+1}{n-4} -3\binom{3n}{n-5}\bigg]
\end{align*}
and
\begin{align*}
[x^{n}]\frac{2t(1 + 6t -9t^{2} + 3t^{3})}{(1-t)^{6}}
& = 2\bigg[\binom{3n+5}{n-1} + 3\binom{3n+4}{n-2} -27\binom{3n+3}{n-3} \\
& \phantom{=}+ 30\binom{3n+2}{n-4}-9\binom{3n+1}{n-5}\bigg],
\end{align*}
with which we find that the variance for the number of returns for S-Motzkin paths of length $3n$ is 
\begin{align*}
\frac{2(313n^{3} + 652n^{2}+53n-178)n}{(2n+5)(2n+4)(2n+3)(2n+2)} & + \frac{n(23n+17)}{2 {(2 n + 3)} {(n + 1)}} - \Big(\frac{n(23n+17)}{2 {(2 n + 3)} {(n + 1)}}\Big)^{2}\\
& = \frac{3  {(14  n^{2} + 31  n + 8)} {(3  n + 2)} {(3  n + 1)} {(n - 1)} n}{4 {(2  n + 5)} {(2  n + 3)}^{2} {(n + 2)} {(n + 1)}^{2}}.
\end{align*}
Similarly, the variance for the number of returns for T-Motzkin paths of length $3n$ is given by 
\begin{align*}
\frac{3 {(79n^{3} + 252 n^{2} + 91 n - 142)} n}{2 {(2 n + 5)} {(2 n + 3)} {(n + 3)} {(n + 2)}} & + \frac{{(19 n + 26)} n}{2(2n+3)(n+2)} - \Big(\frac{{(19 n + 26)} n}{2(2n+3)(n+2)}\Big)^{2}\\
& = \frac{3  {(14  n^{3} + 45  n^{2} + 19  n - 18)}{(3  n + 4)} {(3  n + 2)} n}{4  {(2 n + 5)} {(2  n + 3)}^{2} {(n + 3)}{(n + 2)}^{2}}.
\end{align*}

\subsubsection{Limiting distributions}
\label{sec:lim}

We have defined $t$ implicitly by $t(1-t)^2 = x$. It is well known that this type of implicit equation leads to a square root singularity \cite[Section VII.4]{Flajolet}. In this particular case, the singularity occurs at $x = \frac{4}{27}$, $t = \frac{1}{3}$, where $\frac{d}{dt} t(1-t^2) = (1-t)(1-3t) = 0$. At this point, the singular expansion of $t$ with respect to $x$ is
\begin{align*}
t = \frac13 - \frac{2}{3\sqrt{3}} \Big( 1 - \frac{27x}{4} \Big)^{1/2} + \mathcal{O}\Big( 1 - \frac{27x}{4} \Big).
\end{align*} 
The generating function for the number of returns in S-Motzkin paths is given by
\begin{align*}
S(z, u) & = \frac{(1-t)tu^{2}}{1-tu-tu^{2} + t^{2}u^{2}}.
\end{align*}
Note that for $|x| \leq \frac{4}{27}$ and $|u| \leq 1$, we have $|t| \leq \frac13$ and thus 
\begin{align*}
|1 - tu - tu^2 + t^2 u^2| \geq 1 - |t||u| - |t||u|^2 - |t|^2|u|^2 \geq 1 - \frac13 - \frac13 - \frac19 = \frac29 > 0,
\end{align*}
so the denominator is nonzero and the singularity of $t$ remains the dominant singularity. 
This generating function has the Taylor expansion (with substitution of the singular expansion of $t$):
\begin{align*}
S(z, u) & = \frac{2u^2}{9 - 3u - 2u^2} + \frac{9u^2}{27 - 27u + 4u^3}\Big(t - \frac{1}{3}\Big) + \mathcal{O}\Big(\Big(t - \frac{1}{3}\Big)^{2}\Big) \\
& = \frac{2u^2}{9 - 3u - 2u^2} - \frac{2\sqrt{3}\,u^2}{27 - 27u + 4u^3}\Big(1 - \frac{27x}{4}\Big)^{\frac{1}{2}} + \mathcal{O}\Big(1 - \frac{27x}{4}\Big).
\end{align*}
Applying singularity analysis \cite[Section VI]{Flajolet}, we obtain
\begin{align*}
[x^{n}]S(z, u) & = \frac{2\sqrt{3}\,u^2}{27 - 27u + 4u^3}\cdot \frac{1}{2\sqrt{\pi}}\cdot n^{-3/2}\Big(\frac{27}{4}\Big)^{n}.
\end{align*}
Therefore, the probability generating function for the number of returns in S-Motzkin paths of length $3n$, which is given by $[x^n]S(x,u)/[x^n] S(x,1)$,
converges to
\begin{align*}
\frac{4u^2}{(2u - 3)^2(u + 3)} & = \frac{4}{27}u^2 + \frac{4}{27}u^{3} + \frac{4}{27}u^{4} + \frac{92}{729}u^5 + \cdots. 
\end{align*}
By \cite[Theorem IX.1]{Flajolet}, the distribution of the number of returns in S-Motzkin paths converges to the discrete distribution given by this probability generating function. The probability that the number of returns is precisely $k$ converges to
\begin{align*}
[u^{k}] \frac{4u^2}{(2u - 3)^2(u + 3)} & = \frac{4}{3^{k+3}}(3k\cdot 2^{k-1} - 2^{k} + (-1)^{k}).
\end{align*}

In a similar manner, we find that the limiting probability generating function for the number of returns in T-Motzkin paths of length $3n$ is given by 
\begin{align*}
\frac{4u}{(2u - 3)^2(u + 3)} & = \frac{4}{27}u + \frac{4}{27}u^{2} + \frac{4}{27}u^{3} + \frac{92}{729}u^4 + \frac{76}{729}u^{5} + \cdots,
\end{align*}
and the probability that the number of returns is precisely $k$ converges to 
\begin{align*}
[u^{k}] \frac{4u}{(2u - 3)^2(u + 3)} & = \frac{4}{3^{k+4}}(3k\cdot 2^{k} + 2^{k} - (-1)^{k}).
\end{align*}
The convergence in both cases is demonstrated in the figures below. 

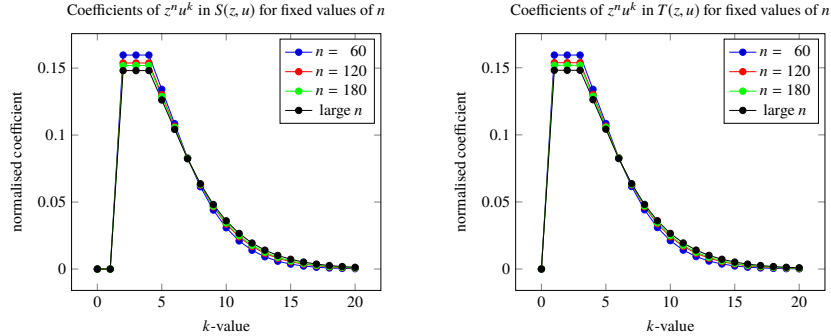
\begin{minipage}{0.5\textwidth}
\begin{tikzpicture}[scale = 0.6]
  \begin{axis}[
    title = {Coefficients of $z^{n}u^{k}$ in $S(z,u)$ for fixed values of $n$},
    xlabel=$k$-value,
    ylabel=normalised coefficient,
    tick label style={/pgf/number format/fixed}]
    \addplot coordinates {(0, 0)  (1, 0)  (2, 27415449/171614024)  (3, 27415449/171614024)  (4, 27415449/171614024)  (5, 84428106/629457851)  (6, 23613271/217621901) (7, 279241093/3369075028)  (8, 4735914/77145841)  (9, 10319007/234403361)  (10, 9565947/310748093)  (11, 8801688/419042245) (12, 29316721/2091632233)  (13, 5787907/632490665)  (14, 17739728/3032877625)  (15, 5445483/1487840839)  (16, 6576494/2933788717) (17, 3214070/2392878923)  (18, 2662909/3384191200)  (19, 3360828/7464015169)  (20, 2987864/11884635417)};
    \addlegendentry{$n= \phantom{1}60$}
    \addplot [red, mark = *] coordinates {(0, 0)  (1, 0)  (2, 44415619/288742649)  (3, 44415619/288742649)  (4, 44415619/288742649)  (5, 77851139/598213173)  (6, 42581871/400000000) (7, 32890610/397374067)  (8, 45675962/728070345)  (9, 31036003/669573146)  (10, 5816075/173001018)  (11, 19134521/798175858) (12, 15162780/899819879)  (13, 14431727/1234562495)  (14, 38469733/4800859044)  (15, 16322231/3005000210)  (16, 6524332/1790774725)  (17, 5974651/2469641898)  (18, 4712001/2961835189) (19, 2204808/2127473537)  (20, 3393865/5073934574)};
    \addlegendentry{$n= 120$}
    \addplot [green, mark = *] coordinates{(0, 0)  (1, 0)  (2, 33872457/222984344)  (3, 33872457/222984344)  (4, 33872457/222984344)  (5, 39686229/308074771)  (6, 86169631/814957291) (7, 9462637/114490218)  (8, 19910334/315661421)  (9, 8379607/178253418) (10, 11314697/328401769)  (11, 3856548/155043103)  (12, 10562709/595541831)  (13, 294986/23586127)  (14, 9843393/1127175191) (15, 9842059/1628659051)  (16, 5713697/1377668843)  (17, 3749507/1327493267)  (18, 3688950/1931697911)  (19, 361389/281825414) (20, 3388318/3961129023)};
    \addlegendentry{$n= 180$}
    \addplot [black, mark = *] coordinates {(0, 0) (1, 0) (2, 4/27) (3, 4/27) (4, 4/27) (5, 92/729) (6, 76/729) (7,
20/243) (8, 1252/19683) (9, 316/6561) (10, 236/6561) (11, 14108/531441) (12, 10316/531441) (13, 7484/531441) (14, 145636/14348907) (15,
104372/14348907) (16, 24812/4782969) (17, 1427228/387420489) (18, 336580/129140163) (19, 237332/129140163) (20, 13514980/10460353203)};
    \addlegendentry{large $n$}
  \end{axis}
\end{tikzpicture}
\end{minipage}
\begin{minipage}{0.5\textwidth}
\begin{tikzpicture}[scale = 0.6]
  \begin{axis}[
    title = {Coefficients of $z^{n}u^{k}$ in $T(z,u)$ for fixed values of $n$},
    xlabel=$k$-value,
    ylabel=normalised coefficient,
    tick label style={/pgf/number format/fixed}]
    \addplot coordinates {(0, 0)  (1, 63384365/397422177)  (2, 63384365/397422177)  (3,
63384365/397422177)  (4, 29282450/218599427)  (5, 30021912/276901465) 
(6, 16560815/199799887)  (7, 274972276/4474121987)  (8,
14313640/324314751)  (9, 12388854/400713899)  (10, 8476223/400937541) 
(11, 12997597/918946744)  (12, 49362281/5329074969)  (13,
8786699/1478753166)  (14, 23446413/6279725401)  (15, 3088268/1343975619)
(16, 2879178/2079481709)  (17, 2970329/3638701445)  (18,
1000181/2125533262)  (19, 655193/2472991671)  (20, 493529/3391286741)};
    \addlegendentry{$n= \phantom{1}60$}
    \addplot [red, mark = *] coordinates {(0, 0)  (1, 54576883/354946946)  (2, 54576883/354946946)  (3,
54576883/354946946)  (4, 23455350/180292669)  (5, 53215633/500000000) 
(6, 35111563/424224022)  (7, 55645731/886818178)  (8,
25350427/546646430)  (9, 15388019/457336432)  (10, 32861068/1369031329) 
(11, 4333164/256696085)  (12, 15030374/1282794435)  (13,
8767853/1090957853)  (14, 8602658/1577976655)  (15, 7095061/1938724016) 
(16, 1316675/541338909)  (17, 2072948/1294755663)  (18,
658451/630657984)  (19, 2217197/3286405807)  (20, 1032407/2389629013)};
    \addlegendentry{$n= 120$}
    \addplot [green, mark = *] coordinates{(0, 0)  (1, 49212484/324028303)  (2, 49212484/324028303)  (3,
49212484/324028303)  (4, 72736823/564723415)  (5, 9809721/92785721)  (6,
26912091/325622447)  (7, 21559927/341789355)  (8, 8168083/173719668) 
(9, 6078049/176350687)  (10, 1602573/64394123)  (11, 9538868/537429265) 
(12, 4494299/359009144)  (13, 16055079/1836267295)  (14,
10318151/1704899510)  (15, 8635829/2078486207)  (16, 5880946/2077639313)
(17, 6533794/3412731633)  (18, 3124474/2429425313)  (19,
2490845/2902079829)  (20, 1584068/2783001489)};
    \addlegendentry{$n= 180$}
    \addplot [black, mark = *] coordinates {(0, 0) (1, 4/27) (2, 4/27) (3, 4/27) (4, 92/729) (5, 76/729) (6, 20/243)
(7, 1252/19683) (8, 316/6561) (9, 236/6561) (10, 14108/531441) (11,
10316/531441) (12, 7484/531441) (13, 145636/14348907) (14,
104372/14348907) (15, 24812/4782969) (16, 1427228/387420489) (17,
336580/129140163) (18, 237332/129140163) (19, 13514980/10460353203) (20,
9476020/10460353203)};
    \addlegendentry{large $n$}
  \end{axis}
\end{tikzpicture}
\end{minipage}

\subsection{The number of peaks}
\label{sec:peaks}

There are two possible types of peaks:
\begin{equation*}
(1) \quad \begin{tikzpicture}[scale = 0.25, line width = 0.2mm]
\coordinate (aux) at (0,0);
\foreach \i in {1, -1}
    \draw[line cap = round] (aux)--++(1,\i) coordinate (aux);
\end{tikzpicture}\qquad\text{and}\qquad
(2) \quad
\begin{tikzpicture}[scale = 0.25, line width = 0.2mm]
\coordinate (aux) at (0,0);
\foreach \i in {1, 0, -1}
    \draw[line cap = round] (aux)--++(1,\i) coordinate (aux);
\end{tikzpicture}
\end{equation*}
We first consider peaks of type (1) and again use the variable $u$ to count them. 
Then from the symbolic equations given below, we obtain the results in Table \ref{peak1}.
\begin{align*}
S(z, u) & = uz^{3} + z^{3}T(z, u)S(z, u)+uz^{3}S(z, u)+z^{3}S(z, u) + z^{3}S(z, u)^{2},\\
T(z, u) & = 1 + z^{3}T(z, u)^{2} + uz^{3}T(z, u) + z^{3}S(z, u)T(z, u)
\end{align*}
\begin{table}[h!]
\caption{Results for peaks of type $(1)$}
\setlength{\tabcolsep}{3mm}
\def\arraystretch{1.5}
\normalsize
\centering
\begin{tabular}{|c||c|c|}
\hline
$K(z, u)$ & $S(z, u)$ & $T(z, u)$\\
\hline
\hline
$\frac{\partial}{\partial u}K(z, u)\Big|_{u=1}$ & $\frac{t ( 1-2t ) }{ ( 1-3t )  ( 1-t ) }$ & $\frac{t}{ ( 1-3t)  ( 1-t ) }$\\
\hline
$[x^{n}]\frac{\partial}{\partial u}K(z, u)\Big|_{u=1}$ & $\binom{3n}{n-1}-2\binom{3n-1}{n-2}$ & $\binom{3n}{n-1}$ \\
\hline
Mean & $\frac{n}{3}+\frac{2}{3}$ & $\frac{n(n+1)}{(3n+1)}$\\
\hline
\hline
$\frac{\partial^{2}}{(\partial u)^{2}}K(z, u)\Big|_{u=1}$ & $\frac{2t^{2} {(1 - 5t + 8t^{2} - 3t^{3})}}{{(1 - 3t)}^{3} {(1 - t)}}$ & $\frac{2t^{2} {(1 - 2t)}}{{(1 - 3t)}^{3} {(1 - t)}}$\\
\hline
$[x^{n}]\frac{\partial^{2}}{(\partial u)^{2}}K(z, u)\Big|_{u=1}$ & $\binom{3n-2}{n-3}\frac{n(n+3)}{(n-2)}$ & $\binom{3n-1}{n-2}n$ \\
\hline
Variance & $\frac{2 {(2 n + 1)} {(n - 1)}}{9 {(3n - 1)}}$ & $\frac{2 {(2n + 1)} {(n + 1)} n}{3 {(3 n + 1)}^{2}}$\\
\hline
\end{tabular}
\label{peak1}
\end{table}

The system of equations for $S(z,u)$ and $T(z,u)$ satisfies the technical conditions of \cite{Drmota}, where it is shown that we have convergence to a normal law in a rather general setting. By the main result of \cite{Drmota}, the number of peaks (of both types) asymptotically follows a Gaussian distribution. 

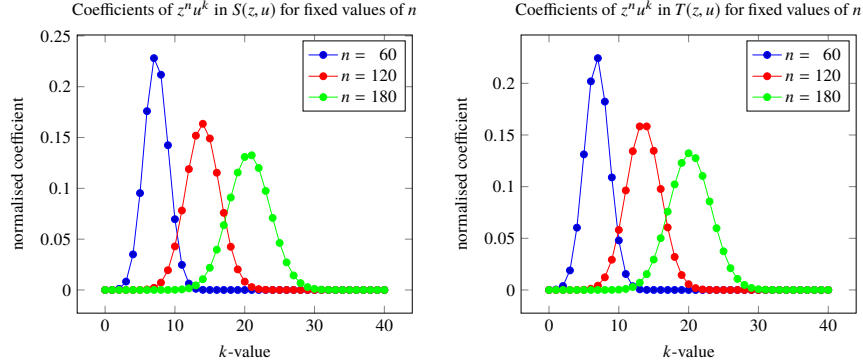
\begin{minipage}{0.5\textwidth}
\begin{tikzpicture}[scale = 0.65]
  \begin{axis}[
    title = {Coefficients of $z^{n}u^{k}$ in $S(z,u)$ for fixed values of $n$},
    xlabel=$k$-value,
    ylabel=normalised coefficient,
    tick label style={/pgf/number format/fixed}]
    \addplot coordinates {(0, 0)  (1, 959835/14950005424)  (2, 5600438/5050168235)  (3,
6021404/730320253)  (4, 22710999/648131137)  (5, 26728558/280435683) 
(6, 54967442/312388747)  (7, 94050499/412331656)  (8,
101286469/478213342)  (9, 2781611/19531250)  (10, 15908217/228478313) 
(11, 10056661/407048629)  (12, 12744618/2017527497)  (13,
4042445/3519649758)  (14, 780357/5362692586)  (15, 307148/24625443049) 
(16, 126581/182674273410)  (17, 12591/537849539599)  (18,
6883/15828130373962)  (19, 443/119190443900833)  (20,
11/1124641208846245)  (21, 0)  (22, 0)  (23, 0)  (24, 0)  (25, 0)  (26,
0)  (27, 0)  (28, 0)  (29, 0)  (30, 0)  (31, 0)  (32, 0)  (33, 0)  (34,
0)  (35, 0)  (36, 0)  (37, 0)  (38, 0)  (39, 0)  (40, 0)};
    \addlegendentry{$n= \phantom{1}60$}
    \addplot [red, mark = *] coordinates {(0, 0)  (1, 1532/826306318149)  (2, 39537/574130624372)  (3,
137554/115912553141)  (4, 239852/18975324097)  (5, 319945/3420498987) 
(6, 2076177/4051686227)  (7, 5212785/2410701229)  (8, 5409626/749181833)
(9, 12951913/665849479)  (10, 21999548/513046519)  (11,
2502187/32000000)  (12, 30977191/260464791)  (13, 24596633/161988334) 
(14, 17850344/109161623)  (15, 31342417/210237165)  (16,
92165373/798933665)  (17, 28220899/371841081)  (18, 19087474/449234209) 
(19, 25941806/1281440863)  (20, 11851131/1444511608)  (21,
5734317/2030268494)  (22, 7879651/9558861027)  (23, 1337843/6577205936) 
(24, 746165/17646458007)  (25, 408085/55351449571)  (26,
94019/87673167561)  (27, 20479/158412216968)  (28, 2261/176433795521) 
(29, 2536/2450828871049)  (30, 94/1397055626497)  (31,
239/68781345052914)  (32, 58/413954337629127)  (33, 9/2084045976231832) 
(34, 1/10097698956385554)  (35, 1/613074579419665165)  (36,
1/54359279380320346991)  (37, 1/7534196118934264549597)  (38,
1/1811974167316313035155926)  (39, 1/906591074861087078765401480)  (40,
1/1414282077065374245500055610172) };
    \addlegendentry{$n= 120$}
    \addplot [green, mark = *] coordinates{(0, 0)  (1, 31/583800870433049)  (2, 202/66625862064237)  (3,
5595/67948848722654)  (4, 42517/29987689431757)  (5, 10163/583865457020)
(6, 126841/780753135744)  (7, 44461/37042723539)  (8,
254944/35326922255)  (9, 537981/14930962663)  (10, 1179520/7776661769) 
(11, 2597197/4767722415)  (12, 3567301/2116922322)  (13,
10304509/2277504577)  (14, 12370309/1165865579)  (15,
15762202/721470255)  (16, 14081003/354952009)  (17, 14984769/235034501) 
(18, 30939332/340106467)  (19, 62934115/544718264)  (20,
24878391/190072172)  (21, 48409190/365338323)  (22, 101750987/847646151)
(23, 15014053/154108459)  (24, 886257/12500000)  (25,
17899019/386638864)  (26, 36091405/1330295068)  (27,
19841505/1390158892)  (28, 13888381/2060606344)  (29,
4330106/1516093753)  (30, 3847585/3544415184)  (31, 2050844/5545857301) 
(32, 1773419/15720910299)  (33, 723111/23493525788)  (34,
269263/35895491830)  (35, 82445/50570000827)  (36, 36515/115777952304) 
(37, 45898/846976392295)  (38, 9464/1147204910935)  (39,
9125/8223702924342)  (40, 777/5911234528694) };
    \addlegendentry{$n= 180$}
  \end{axis}
\end{tikzpicture}
\end{minipage}
\begin{minipage}{0.5\textwidth}
\begin{tikzpicture}[scale = 0.65]
  \begin{axis}[
    title = {Coefficients of $z^{n}u^{k}$ in $T(z,u)$ for fixed values of $n$},
    xlabel=$k$-value,
    ylabel=normalised coefficient,
    tick label style={/pgf/number format/fixed}]
    \addplot coordinates {(0, 383371/17345010459)  (1, 3683283/8332213723)  (2,
3916570/1025888793)  (3, 12063790/637530763)  (4, 25141599/416831066) 
(5, 11010833/83893530)  (6, 30768806/152381517)  (7, 75718036/337492013)
(8, 82548069/452843051)  (9, 109303370/1003206299)  (10,
27411841/571797465)  (11, 10712527/692719870)  (12, 8231863/2271185278) 
(13, 2284763/3755955025)  (14, 1426398/19931449019)  (15,
119859/20935281982)  (16, 66319/222406442056)  (17, 6887/726373725338) 
(18, 603/3625118530705)  (19, 157/116565896723502)  (20,
11/3266814940506806)  (21, 0)  (22, 0)  (23, 0)  (24, 0)  (25, 0)  (26,
0)  (27, 0)  (28, 0)  (29, 0)  (30, 0)  (31, 0)  (32, 0)  (33, 0)  (34,
0)  (35, 0)  (36, 0)  (37, 0)  (38, 0)  (39, 0)  (40, 0)};
    \addlegendentry{$n= \phantom{1}60$}
    \addplot [red, mark = *] coordinates {(0, 2153/3427108480803)  (1, 15767/627441005719)  (2,
46861/100413029533)  (3, 298830/55737080891)  (4, 644630/15050741093) 
(5, 923528/3642300155)  (6, 7629356/6591011701)  (7, 2098038/501102407) 
(8, 23049060/1884103363)  (9, 23939880/817237237)  (10,
170538406/2934314833)  (11, 24181487/250984803)  (12, 10985317/81779105)
(13, 27232647/172021531)  (14, 27232647/172021531)  (15,
48016869/356454110)  (16, 48935575/500759431)  (17, 30620483/506044845) 
(18, 23821653/744578143)  (19, 10605691/734399639)  (20,
1722315/309774268)  (21, 26925999/14770813688)  (22, 2755367/5425538145)
(23, 1797191/14993391562)  (24, 711620/29800502171)  (25,
97099/24292636717)  (26, 57407/102690623837)  (27, 24657/379949123954) 
(28, 8987/1448756875745)  (29, 1418/2932104803973)  (30,
491/16152141642100)  (31, 175/115189801101856)  (32, 21/353863068863782)
(33, 5/2818968693980943)  (34, 1/25330447254111666)  (35,
1/1583152953582489206)  (36, 1/144383549384643785828)  (37,
1/20567436608745240960534)  (38, 1/5080156843258651931261057)  (39,
1/2608660538922538977831486594)  (40, 1/4173856862058298330836381741131)};
    \addlegendentry{$n= 120$}
    \addplot [green, mark = *] coordinates{(0, 14/782311689950087)  (1, 84/78231168980437)  (2,
4263/139070267443888)  (3, 6551/11803433442755)  (4, 6828/952645640881) 
(5, 31281/444364113772)  (6, 85578/156302150173)  (7,
292074/84238788379)  (8, 347949/19074979334)  (9, 848294/10478696593) 
(10, 3202495/10441763466)  (11, 4792891/4786232556)  (12,
30723755/10819774149)  (13, 8489593/1206312054)  (14, 4012691/261835890)
(15, 35335609/1199783461)  (16, 38942637/776745733)  (17,
14811603/195312500)  (18, 49460009/483980811)  (19, 60925570/495491831) 
(20, 45372754/342861621)  (21, 22419783/175717738)  (22,
32683905/296230706)  (23, 34165611/398880725)  (24, 49806429/833766556) 
(25, 28384041/758032805)  (26, 17032029/807200491)  (27,
12697003/1187824577)  (28, 1945172/399628949)  (29, 2068987/1038916825) 
(30, 3756592/5134156929)  (31, 1839757/7627034509)  (32,
2744995/38508389698)  (33, 409571/21716201274)  (34, 211877/47492327049)
(35, 172299/182926628125)  (36, 52621/297039117900)  (37,
16533/558249603341)  (38, 7046/1605055266029)  (39, 234/406735888387) 
(40, 171/2573420539400)};
    \addlegendentry{$n= 180$}
  \end{axis}
\end{tikzpicture}
\end{minipage}

We now consider peaks of type $(2)$, and again use the variable $u$ to count them. From the symbolic equations we obtain
\begin{align*}
S(z, u) & = z^{3}+z^{3}T(z, u)S(z, u)+uz^{3}S(z, u)+z^{3}S(z, u)+z^{3}S(z, u)^{2},\\
T(z, u) & = 1 + z^{3}T(z, u)^{2}+uz^{3}T(z, u)+z^{3}S(z, u)T(z, u).
\end{align*}
Note that $\mathcal{T}$ contains an empty path. As a result, for paths of the form 
\begin{align*}
\begin{tikzpicture}[scale=0.4]
\draw [line cap = round](0,0) to  (0.5,0.5);
\draw [line cap = round](0.5,0.5) to  (1,0.5);
	\node at (1.6, 0.8) {$\mathcal{T}$};	
\draw [line cap = round](1.8,0.5) to  (2.3,0);
	\node at (2.8, 0.3) {$\mathcal{U}$};
\end{tikzpicture}\qquad \raisebox{1.5ex}{\text{and}} \qquad
\begin{tikzpicture}[scale=0.4]
\draw [line cap = round](0,0) to  (0.5,0.5);
\draw [line cap = round](0.5,0.5) to  (1,0.5);
	\node at (1.6, 0.8) {$\mathcal{T}$};	
\draw [line cap = round](1.8,0.5) to  (2.3,0);
	\node at (2.8, 0.3) {$\mathcal{T}$};
\end{tikzpicture}\, ,
\end{align*}
if the T-Motzkin path is empty we obtain $uz^{3}S(z, u)$ and $uz^{3}T(z, u)$ respectively. If path is not empty we obtain  $z^{3}T(z, u)S(z, u)$ and $z^{3}T(z, u)^{2}$. 
Using the generating function equations we obtain the following results. 
\begin{table}[h!]
\caption{Results for peaks of type $(2)$}
\setlength{\tabcolsep}{3mm}
\def\arraystretch{1.75}
\normalsize
\centering
\begin{tabular}{|c||c|c|}
\hline
$K(z, u)$ & $S(z, u)$ & $T(z, u)$\\
\hline
\hline
$\frac{\partial}{\partial u}K(z, u)\Big|_{u=1}$ & $\frac{t^{2}}{1-3t}$ & $\frac{t}{1-3t}$ \\
\hline
$[x^{n}]\frac{\partial}{\partial u}K(z, u)\Big|_{u=1}$ & $\binom{3n-2}{n-2}$ & $\binom{3n-1}{n-1}$\\
\hline
Mean & $\frac{(2n+1)(n-1)}{3(3n-1)}$ & $\frac{(2n+1)(n+1)}{3(3n+1)}$ \\
\hline
\hline
$\frac{\partial^{2}}{(\partial u)^{2}}K(z, u)\Big|_{u=1}$ & $\frac{2  {(1-2t)} {(1-t)}t^{3}}{{(1-3t)}^{3}}$ & $\frac{2{(1 - 3t + 3t^{2})} {(1-t)} t^{2}}{{(1-3t)}^{3}}$ \\
\hline
$[x^{n}]\frac{\partial^{2}}{(\partial u)^{2}}K(z, u)\Big|_{u=1}$ & $\binom{3n-3}{n-3}\frac{2n}{3}$ & $\binom{3n-3}{n-2}n$\\
\hline
Variance & $\frac{2 {(10 n^{2} - 11 n + 2)} {(2n + 1)} {(n - 1)}}{9 {(3n - 1)}^{2}{(3 n - 2)}}$ & $\frac{2 {(30 n^{3} - 23 n^{2} - 3 n + 2)}{(2n + 1)} {(n + 1)}}{9 {(3 n + 1)}^{2} {(3 n - 1)} {(3 n - 2)}}$ \\
\hline
\end{tabular}
\end{table}

\subsection{Valleys}
\label{sec:valleys}
There are two possible types of valleys:
\begin{equation*}
(1) \quad \begin{tikzpicture}[scale = 0.25, line width = 0.2mm]
\coordinate (aux) at (0,0);
\foreach \i in {-1, 1}
    \draw[line cap = round] (aux)--++(1,\i) coordinate (aux);
\end{tikzpicture}\qquad\text{and}\qquad
(2) \quad
\begin{tikzpicture}[scale = 0.25, line width = 0.2mm]
\coordinate (aux) at (0,0);
\foreach \i in {-1, 0, 1}
    \draw[line cap = round] (aux)--++(1,\i) coordinate (aux);
\end{tikzpicture}
\end{equation*}
For valleys of type $(1)$, using the variable $u$ to count them we obtain generating function equations 
\begin{align*}
S(z, u) & = z^{3} + uz^{3}T(z, u)S(z, u)+z^{3}S(z, u)+z^{3}S(z, u)+z^{3}S(z, u)^{2},\\
T(z, u) & = 1 + uz^{3}T(z, u)(T(z, u)-1) + 2z^{3}T(z, u)+z^{3}S(z, u)T(z, u).
\end{align*}
This is taking into account that the empty T-Motzkin path does not contribute a valley of type $(1)$. 
From the generating function equations we obtain: 
\begin{table}[h!]
\caption{Results for valleys of type $(1)$}
\setlength{\tabcolsep}{3mm}
\def\arraystretch{1.75}
\normalsize
\centering
\begin{tabular}{|c||c|c|}
\hline
$K(z, u)$ & $S(z, u)$ & $T(z, u)$\\
\hline
\hline
$\frac{\partial}{\partial u}K(z, u)\Big|_{u=1}$ & $\frac{t^{2}}{{(1-3t)} {(1-t)}}$ & $\frac{2t^{2}}{{(1-3t)} {(1-t)}^{2}}$\\
\hline
$[x^{n}]\frac{\partial}{\partial u}K(z, u)\Big|_{u=1}$ & $\binom{3n-1}{n-2}$ & $2\binom{3n}{n-2}$ \\
\hline
Mean & $\frac{n}{3}-\frac{1}{3}$ & $\frac{n(n-1)}{3n+1}$ \\
\hline
\hline
$\frac{\partial^{2}}{(\partial u)^{2}}K(z, u)\Big|_{u=1}$ & $\frac{2(1 - t - 3t^{2})t^3}{(1-3t)^{3}(1-t)}$ & $\frac{2(2 - t - 9t^{2})t^{3}}{(1-3t)^{3}(1-t)^{2}}$ \\
\hline
$[x^{n}]\frac{\partial^{2}}{(\partial u)^{2}}K(z, u)\Big|_{u=1}$ & $(n-1)\binom{3n-2}{n-3}$ & $\frac{(n-1)(n-2)}{(n+1)}\binom{3n-1}{n-2}$ \\
\hline
Variance & $\frac{2{(2 n + 1)} {(n - 1)}}{9{(3n - 1)}}$ & $\frac{2 {(2 n + 1)} {(n + 1)} {(n - 1)}}{3 {(3 n + 1)}^{2}}$ \\
\hline
\end{tabular}
\end{table}

As in our analysis of peaks in the previous subsection, we can apply the main result of \cite{Drmota} to prove that the number of valleys (of both types) asymptotically follows a Gaussian distribution. 

The generating function equations for valleys of type $(2)$, again using $u$ to count the number of valleys, are given by 
\begin{align*}
S(z, u) & = z^{3} + z^{3}T(z, u)S(z, u)+uz^{3}S(z, u)+z^{3}S(z, u)+uz^{3}S(z, u)^{2},\\
T(z, u) & = 1 + z^{3}T(z, u)^{2} + uz^{3}(T(z, u)-1) + z^{3} + uz^{3}S(z, u)(T(z, u)-1)\\
&\phantom{=} + z^{3}S(z, u).
\end{align*}
Again, we take into account the absence of a valley in the case of an empty T-Motzkin path. 
The equations yield:
\begin{table}[h!]
\caption{Results for valleys of type $(2)$}
\setlength{\tabcolsep}{3mm}
\def\arraystretch{1.75}
\normalsize
\centering
\begin{tabular}{|c||c|c|}
\hline
$K(z, u)$ & $S(z, u)$ & $T(z, u)$\\
\hline
\hline
$\frac{\partial}{\partial u}K(z, u)\Big|_{u=1}$ & $\frac{t^2}{(1-3t)}$ & $\frac{2t^2}{(1-3t)(1-t)}$ \\	
\hline
$[x^{n}]\frac{\partial}{\partial u}K(z, u)\Big|_{u=1}$ & $\binom{3n-2}{n-2}$ & $2\binom{3n-1}{n-2}$\\
\hline
Mean & $\frac{(n-1)(2n+1)}{3(3n-1)}$ & $\frac{2(n+1)(n-1)}{3(3n+1)}$\\
\hline
\hline
$\frac{\partial^{2}}{(\partial u)^{2}}K(z, u)\Big|_{u=1}$ & $\frac{2(1-2t)(1-t)t^3}{(1-3t)^3}$ & $\frac{2(2 - 3t - 3t^2)t^3}{(1-3t)^3}$ \\
\hline
$[x^{n}]\frac{\partial^{2}}{(\partial u)^{2}}K(z, u)\Big|_{u=1}$ & $\frac{2n}{3}\binom{3n-3}{n-3}$ & $2(n-1)\binom{3n-3}{n-3}$ \\
\hline
Variance & $\frac{2 {(10 n^{2} - 11 n + 2)} {(2 n +1)} {(n - 1)}}{9 {(3 n - 1)}^{2}{(3n - 2)}}$ & $\frac{4 {(15 n^{2} - 19n+ 8)} {(2 n +1)} {(n + 1)} {(n - 1)}}{9 {(3 n+ 1)}^{2} {(3n - 1)} {(3n - 2)}}$ \\
\hline
\end{tabular}
\end{table}

\subsection{Valleys on the $x$-axis}

We now consider valleys that lie on the $x$-axis. Keeping the two types of valleys discussed in the previous subsection, a valley of type $(1)$ contributes one return, and a valley of type $(2)$ contributes two returns.  
For valleys on the $x$-axis of type $(1)$, using the variable $u$ to count them we obtain generating function equations 
\begin{align*}
S(z, u) & = z^{3} + uz^{3}T(z, 1)S(z, u)+z^{3}S(z, u)+z^{3}S(z, 1)+z^{3}S(z, u)S(z, 1),\\
T(z, u) & = 1 + uz^{3}T(z, 1)(T(z, u)-1) + z^{3}T(z, 1) + z^{3}T(z, u)+z^{3}S(z, 1)T(z, u).
\end{align*}
Let $v_{1} = 30 n^{3} + 43 n^{2} + 154 n + 288$ and $v_{2} = 778 n^{6} + 3953 n^{5} + 11212 n^{4} + 24373 n^{3} + 30064 n^{2} + 16260 n + 2160$, then we obtain the following results.
\begin{table}[h!]
\caption{Results for valleys on the $x$-axis of type $(1)$}
\setlength{\tabcolsep}{3mm}
\def\arraystretch{1.75}
\normalsize
\centering
\begin{tabular}{|c||c|c|}
\hline
$K(z, u)$ & $S(z, u)$ & $T(z, u)$\\
\hline
\hline
$\frac{\partial}{\partial u}K(z, u)\Big|_{u=1}$ & $\frac{t^2}{(1-t)^3}$ & $\frac{(2-t)t^2}{(1-t)^4}$\\
\hline
$[x^{n}]\frac{\partial}{\partial u}K(z, u)\Big|_{u=1}$ & $\binom{3n+1}{n-2}-3\binom{3n}{n-3}$ & $2\binom{3n+2}{n-2}-7\binom{3n+1}{n-3}+3\binom{3n}{n-4}$\\
\hline
Mean & $\frac{7{(n - 1)} n}{2 {(2n + 3)}{(n + 1)}}$ & $\frac{{(19 n + 18)} {(n - 1)} n}{2{(3n + 1)} {(2 n + 3)} {(n + 2)}}$\\
\hline
\hline
$\frac{\partial^{2}}{(\partial u)^{2}}K(z, u)\Big|_{u=1}$ & $\frac{2t^3}{(1-t)^5}$ & $\frac{2(2-t)t^3}{(1-t)^6}$\\
\hline
$[x^{n}]\frac{\partial^{2}}{(\partial u)^{2}}K(z, u)\Big|_{u=1}$ & $2\binom{3n+2}{n-3}-6\binom{3n+1}{n-4}$ & $4\binom{3n+3}{n-3}-14\binom{3n+2}{n-4}+6\binom{3n+1}{n-5}$\\
\hline
Variance & $ \frac{{v_{1}} {(3n + 1)} {(n - 1)} n}{4 {(2 n + 5)}{(2n + 3)}^{2} {(n + 2)} {(n + 1)}^{2}}
$ & $\frac{{v_{2}} {(n - 1)}n}{4 {(3 n + 1)}^{2} {(2 n + 5)} {(2n + 3)}^{2} {(n + 3)} {(n + 2)}^{2}}$  \\
\hline
\end{tabular}
\end{table}

In a similar manner to that of Section \ref{sec:lim}, we find the limiting probability generating function for valleys of type $(1)$ on the $x$-axis is given by $\frac{4(u+3)}{(7-3u)^{2}}$ for S-Motzkin paths, and $\frac{4(u+11)}{3(7-3u)}$ for T-Motzkin paths (both of length $3n$).  

For valleys on the $x$-axis of type $(2)$, again using $u$ to count the number of valleys, we obtain the generating function equations
\begin{align*}
S(z, u) & = z^{3} + z^{3}T(z, 1)S(z, u)+uz^{3}S(z, u)+z^{3}S(z, 1)+uz^{3}S(z, u)S(z, 1),\\
T(z, u) & = 1 + z^{3}T(z, u)T(z, 1) + z^{3} + uz^{3}(T(z, u)-1) + uz^{3}S(z, 1)(T(z, u)-1)\\
&\phantom{=} + z^{3}S(z, 1).
\end{align*}
From these the following results are obtained. 
\begin{table}[h!]
\caption{Results for valleys on the $x$-axis of type $(2)$}
\setlength{\tabcolsep}{3mm}
\def\arraystretch{1.75}
\normalsize
\centering
\begin{tabular}{|c||c|c|}

\hline
$K(z, u)$ & $S(z, u)$ & $T(z, u)$\\
\hline
\hline
$\frac{\partial}{\partial u}K(z, u)\Big|_{u=1}$ & $\frac{t^2}{(1-t)^2}$ & $\frac{(2-t)t^2}{(1-t)^3}$ \\
\hline
$[x^{n}]\frac{\partial}{\partial u}K(z, u)\Big|_{u=1}$ & $\binom{3n}{n-2}-3\binom{3n-1}{n-3}$ & $2\binom{3n+1}{n-2}-7\binom{3n}{n-3}+3\binom{3n-1}{n-4}$\\
\hline
Mean & $\frac{n-1}{n+1}$ & $\frac{{(11n + 6)}{(n - 1)}}{2 {(3 n + 1)} {(2 n + 3)}}$\\
\hline
\hline
$\frac{\partial^{2}}{(\partial u)^{2}}K(z, u)\Big|_{u=1}$ & $\frac{2t^3}{(1-t)^3}$ & $\frac{2(2-t)t^3}{(1-t)^4}$\\
\hline
$[x^{n}]\frac{\partial^{2}}{(\partial u)^{2}}K(z, u)\Big|_{u=1}$ & $2\binom{3n}{n-3}-6\binom{3n-1}{n-4}$ & $4\binom{3n+1}{n-3}-14\binom{3n}{n-4}+6\binom{3n-1}{n-5}$\\
\hline
Variance & $\frac{{(3 n + 1)} {(n - 1)} n}{{(2 n + 3)} {(n + 1)}^{2}}$ & $\frac{{(203n^{3} + 437n^{2} + 268n + 12)}{(n - 1)} n}{4 {(3n + 1)}^{2} {(2n+ 3)}^{2} {(n + 2)}}$\\
\hline
\end{tabular}
\end{table}

The limiting probability generating functions for the number of valleys of type $(2)$ on the $x$-axis is given by $\frac{4}{(3-u)^{2}}$ for S-Motzkin paths and $\frac{13-u}{3(3-u)^{2}}$ for T-Motzkin paths (both of length $3n$).

\section{Identities}
\label{sec:6}

In Sections \ref{sec:peaks} and \ref{sec:valleys} the coefficients of the generating functions used to find the variance were greatly simplified by using derivatives (compared to extracting coefficients using Cauchy's integral formula). An example of this simplification is given: 	
Using Cauchy's integral formula as in Section \ref{sec:5.1.1} we obtain the coefficients
\begin{align*}
[x^{n}]\frac{2t^{2} {\left(1 - 2 t\right)}}{{\left(1 - 3 t\right)}^{3} {\left(1 - t\right)}}& = 2\sum_{k\geq 0}(k+1)3^{k}\bigg[ \binom{3n-k-1}{n-k-2} -2\binom{3n-k-2}{n-k-3}\bigg].
\end{align*} 
On the other hand, the generating function can be expressed as a derivative, and by using the formula $\frac{t^{3}}{1-3t} = \sum_{n\geq 3}\binom{3n-3}{n-3}x^{n}$ we find that
\begin{align*}
[x^{n}]\frac{2t^{2} {\left(1 - 2 t\right)}}{{\left(1 - 3 t\right)}^{3} {\left(1 - t\right)}}&  = [x^{n}]\frac{2}{3}\cdot\frac{d}{dx}\frac{t^{3}}{1-3t} = \binom{3n-1}{n-2}n.
\end{align*}
It follows that 
\begin{align*}
2\sum_{k\geq 0}(k+1)3^{k}\bigg[ \binom{3n-k-1}{n-k-2} -2\binom{3n-k-2}{n-k-3}\bigg] & = \binom{3n-1}{n-2}n,
\end{align*}
which is a special case of the more general identity
\begin{align*}
2\sum_{k\ge j}3^k(k+1)\bigg[\binom{3n-k-1}{n-k-2}-2\binom{3n-k-2}{n-k-3}\bigg]
=\binom{3n-j-1}{n-j-2}(n+j)3^{j}.
\end{align*}
This and other beautiful identities such as 
\begin{equation*}
2\sum_{k\ge0}3^k(k+2i)\binom{3n-k+i-4}{n-k-i-1}=\binom{3n+i-3}{n-i}(n-i)
\end{equation*}
and 
\begin{equation*}
2\sum_{k\ge0}3^k(k+2i+1)\binom{3n-k+i-2}{n-k-i-1}=\binom{3n+i-1}{n-i}(n-i).
\end{equation*}
can be proved directly by induction. 
A table of the simplifications used to calculate variances in Section \ref{sec:5} is given in Table \ref{table:deriv}. 
\begin{table}
\normalsize
\caption{Generating functions and their coefficients}
\setlength{\tabcolsep}{3mm}
\def\arraystretch{1.5}
\centering
\begin{tabular}{ |c|c|c| }
\hline
Generating function & In terms of derivatives & Power series expansion\\
\hline
\hline
$\frac{2t^2(1-5t+8t^2-3t^3)}{(1-3t)^3(1-t)}$ & $-\frac23\frac d{dx}\frac{t^3}{1-3t}+2x\frac d{dx}\frac{t^2}{1-3t}$ & $\sum\limits_{n\ge3}\binom{3n-2}{n-3}\frac{n(n+3)}{(n-2)}x^n$ \\
\hline
$\frac{2t^2(1-2t)}{(1-3t)^3(1-t)}$ & $\frac23\frac d{dx}\frac{t^3}{1-3t}$ & $\sum\limits_{n\ge2}\binom{3n-1}{n-2}nx^n$ \\
\hline
$\frac{2t^3(1-2t)(1-t)}{(1-3t)^3}$ & $\frac{2}{3}x\frac d{dx}\frac{t^3}{1-3t}$ & $\sum\limits_{n\ge3}\binom{3n-3}{n-3}\frac{2n}{3}x^n$ \\
\hline
$\frac{2t^2(1-3t+3t^2)(1-t)}{(1-3t)^3}$ & $x\frac d{dx}\frac{t^2}{1-3t}-x\frac d{dx}\frac{t^3}{1-3t}$ & $\sum\limits_{n\ge2}\binom{3n-3}{n-2}nx^n$ \\
\hline
$\frac{2t^3(1 - t - 3t^{2})}{(1-3t)^{3}(1-t)}$ & $\frac12\frac d{dx}\frac{t^5}{1-3t}+\frac{1}{2}\frac d{dx}\frac{t^4}{1-3t}$ & $\sum\limits_{n\ge3}\binom{3n-2}{n-3}(n-1)x^n$ \\
\hline
$\frac{2t^{3}(2 - t - 9t^{2})}{(1-3t)^{3}(1-t)^{2}}$ & $\frac1x\frac d{dx}\Big(\frac{4}{5}\frac{t^5}{1-3t}+\frac{3}{5}\frac{t^6}{1-3t}-\frac{t^7}{1-3t}\Big) $ & $\sum\limits_{n\ge2}\binom{3n-1}{n-2}\frac{(n-1)(n-2)}{n+1}x^n$ \\
\hline
$\frac{2t^3(2 - 3t - 3t^2)}{(1-3t)^3}$ & $-\frac25\frac d{dx}\frac{t^6}{1-3t}-\frac15\frac d{dx}\frac{t^5}{1-3t}+\frac d{dx}\frac{t^4}{1-3t}$ & $\sum\limits_{n\ge3}2\binom{3n-3}{n-3}(n-1)x^n$ \\
 \hline
\end{tabular}
\label{table:deriv}
\end{table}

\newpage
\bibliography{paper}

\end{document}